\documentclass[12pt]{article}
\usepackage{amsmath}
\usepackage{amssymb}

\usepackage[german,english]{babel}
\usepackage{amsmath, latexsym}
\usepackage{amssymb}
\usepackage{url}
\usepackage{graphicx}
\usepackage[usenames]{color}
\usepackage{gastex}
\usepackage{longtable}
\usepackage{lscape}
\usepackage{tabularx}
\usepackage{multicol}
\usepackage{amsmath}
\usepackage{latexsym}
\usepackage{verbatim}
\usepackage{multirow}
\usepackage{ragged2e}
\usepackage{gastex}
\usepackage{graphicx}

\newtheorem{theo}{{\sc Theorem}}[section]

\newtheorem{remark}{{\em Remark}}

\newtheorem{example}{{\sc Example}}
\begin{document}
\title{ Some Centrally Symmetric  Manifolds }
\author{{Dipendu Maity  and Ashish Kumar Upadhyay}\\[2mm]
{\normalsize Department of Mathematics}\\{\normalsize Indian Institute of Technology Patna }\\{\normalsize Patliputra Colony, Patna 800\,013,  India.}\\
{\small \{dipendumaity, upadhyay\}@iitp.ac.in}}
\maketitle

\vspace{-5mm}

\hrule

\begin{abstract}  We show existence of centrally symmetric maps on surfaces all of whose faces are quadrangles and pentagons for each orientable genus $g \geq 0$. We also show existence of centrally symmetric maps on surfaces all of whose faces are hexagons for each orientable genus $g = 2k-1$, $k\in \mathbb{N}$. We enumerate centrally symmetric triangulated manifolds of dimensions 2 and 3 with few vertices.
\end{abstract}

{\small

{\bf AMS classification\,: 52B70, 57M20, 52A25}

{\bf Keywords\,:} Centrally symmetric manifolds,  Maps, Triangulations
}

\bigskip

\hrule

\section{Introduction and definitions}

In this article we consider simplicial complex which are finite and abstract. A simplicial complex {\em K} is called a {\em combinatorial d-manifold} if $lk_{K}(u)$ is a combinatorial {\em (d-1)}-sphere for each vertex {\em u} in {\em K}. A simplicial complex {\em M} is said to be a {\em triangulated d-manifold} if it's geometric carrier $|M|$ is a topological {\em d}-manifold. Let {\em X} and {\em Y} be two simplicial complexes. An {\em isomorphism} from {\em X} to {\em Y} is a bijection $\phi : V(X)\rightarrow V(Y)$ such that for $\sigma \subseteq V(X)$, $\sigma$ is a simplex of {\em X} if and only if $\phi(\sigma)$ is a simplex of {\em Y}. Two simplicial complexes {\em X, Y} are called (simplicially) isomorphic when such an isomorphism exists. An isomorphism from a simplicial complex {\em X} to itself is called an {\em automorphism} of {\em X}. All the automorphisms of {\em X} form a group, which is denoted by {\em Aut(X)}. Let {\em M} be a {\em d}-dimensional combinatorial manifold with the face vector $(f_{0}, f_{1}, \dots, f_{d})$ where $f_{i}$ denotes the number of {\em i}-dimensional faces of {\em M} for $0\leq i \leq d$. Then the number $\chi(M):= \sum_{i=0}^{d}(-1)^{i}f_{i}$ is called {\em Euler characteristic} of {\em M}. When $d = 2$ then  2-manifold is called surface, that is, a surface {\em S} is a connected, compact, 2-manifold without boundary. Let {\em G := (V, E)} be a finite simple graph with vertex set {\em V} and edge set {\em E}. Let {\em G} be a graph all of whose vertices have degree $\geq 3$. A {\em map M} is an embedding of a graph {\em G} on a surface {\em S} such that the closure of components of {\em S $\setminus$ G}, called the faces of {\em M}, are closed 2-cells, that is, each homeomorphic to 2-disk. A map {\em M} is said to be a polyhedral map if the intersection of any two distinct faces is either empty, a common vertex, or a common edge. A polyhedral map {\em M} is called triangulation if each face of {\em M} is a triangle. A group action is a description of symmetries of objects using groups. An {\em involution} is a function, when we applied it twice, it brings one back to the starting point. In this article, we use the involution which is represented by the permutation $(1, 2m)(2, 2m-1) \dots (m, m+1)$ on $2m$ vertices. We call a polyhedral map {\em M} to be {\em centrally symmetric} in short {\em CS} if it is invariant under an involution of its vertex set which fixes no face of the object. Centrally symmetric has wide applications in mathematics, see \cite{matousek}, \cite{steinlein}, \cite{walkup}.

A subset $C\subseteq R^{m}$ is called convex if for each pair of points {\em a, b $\in$ C} the arc  $ta+(1-t)b \subseteq C$, where $t \in [0,1]$. Let {\em A} be a set in $R^{m}$. The smallest convex set containing {\em A} is called the convex hull of {\em A}. A {\em polytope} is a convex hull of a finite set {\em A}. In particular, if a $(n-1)$-polytope contains exactly {\em n} vertices then we call it $(n-1)$-simplex and we denote it by $\triangle^{n-1}$. We recall that a {\em d}-dimensional polytope {\em P $\subset \mathbb{R}^{d}$} is called centrally symmetric if it is invariant under a point reflection through its center (\cite{wiki}). If we consider the center of the object {\em P} is origin then {\em P} is said to be centrally symmetric if {\em P = -P}. So, if {\em d $>$ 0} then the involution {\em I : x $\rightarrow$ -x} of $\mathbb{R}^{d}$ does not fix any non-trivial face of {\em P} of the polytope {\em P} centered at origin and hence, {\em P} has an even number of vertices $n = 2m$. In this context, Gr$\ddot{u}$nbaum observed there is a unique  4-dimensional CS polytope on  10 vertices (see \cite{grunbaum}). This object is also nearly neighborly. He also showed that there is no nearly neighborly centrally symmetric 4-polytopes with $n \geq 12$ vertices. See \cite{lutz0}, for more extensive literature on general properties of centrally symmetric polytopes and some other results related to centrally symmetric polytopes.

Let {\em M} be a triangulation with {\em n} vertices which can be always regarded as a subcomplex of $\triangle^{n-1}$. In general, any polyhedral map {\em M} can be regarded as a subcomplex of some $ d$-polytope, see \cite{kuhnel}. Let {\em P} be a polytope. The boundary of {\em P} is a sphere. We consider sphere which is boundary of a polytope, that is, a sphere is defined as a simplicial complex $\partial(P)$. A simplicial $(d-1)$-sphere {\em S} is called {\em l-neighbourly} if every set of {\em l} (or less) vertices forms a face of {\em S}. The standard dihedral and cyclic group acts on the set $\{1, 2, \dots, 2m\}$ with generators $a_{2m} = (123\dots 2m)$ and  $b_{2m} = (1, 2m)(2, 2m-1) \dots (m, m+1)$ of $D_{2m} = < a_{2m}, b_{2m}>$ and $Z_{2m} = < a_{2m}>$ respectively. Lassmann and Sparla (in \cite{lassmann}) and Lutz (in \cite{lutz0}) have studied  centrally symmetric spheres and product of spheres under standard dihedral and cyclic group actions. Lassmann and Sparla \cite{lassmann} showed that there are three centrally symmetric 3-neighbourly triangulations of the product $S^{2} \times S^{2}$ with cyclic symmetry. Lutz \cite{lutz0} has extended this result and enumerated triangulations of product of spheres using cyclic and dihedral group action on $n = 2d + 4$ vertex where $d$ is the dimension of the sphere.

A triangulation {\em M} on $ n$ vertices (see \cite{kuhnel}) is said to be tight or $ 2$-neighbourly if its edge graph is a complete graph $K_{n}$ and it satisfies $(n-3)(n-4) = 6(2-\chi(M))=12g$ where $g$ is the genus of {\em M}. A tight or {\em 2}-neighbourly triangulation {\em M} on {\em n} vertices is said to be centrally symmetric if it avoids fixing a face under an involution and contains $\binom{n}{2}-\frac{n}{2}$ edges. We denote it by $n_{tight}$. Also, the map {\em M} satisfies the equality $2(\frac{n}{2}-1)(\frac{n}{2}-3) = 3(2-\chi(M))$ (see \cite{kuhnel}). When $n$ = 10 the above  equation is not satisfied. Therefore, there does not exist tight triangulated centrally symmetric surface on $ 10$ vertices. Lutz \cite{lutz0} showed existence of vertex transitive centrally symmetric triangulation of spheres and torus under cyclic and dihedral group action. In this article we relaxed the condition of vertex transitivity. We have extended this result to centrally symmetric triangulation of surfaces for few vertices under $\mathbb{Z}_{2}$ group action. We have enumerated these objects by using computer. In {\em Section 2}, we use an idea which has been introduced in \cite{kuhnel} to construct CS manifold from already known CS manifolds. For each orientable genus we show the existence of centrally symmetric quadrangulations in {\em Section 3}. Again, for each orientable genus we show existence of CS maps all of whose faces are pentagons in {\em Section 4}. Also, for each positive odd orientable genus we show the existence of CS surfaces all of whose faces are hexagons in {\em Scetion 5}. In {\em Section 6} we give an idea to construct centrally symmetric maps on surfaces of type $\{q, p\}$ from known centrally symmetric maps on surfaces of type $\{p, q\}$. In {\em Section 7} we enumerate centrally symmetric triangulated  3-manifolds on  12 vertices by using computer. We use the notation {\em CS} in place of centrally symmetric and {\em CST} in place of centrally symmetric triangulated map throughout this article. The main results of this article are :

\begin{theo}\label{thm1}
There are exactly {\em 6303} centrally symmetric triangulated surfaces with $n \leq 12$ vertices up to isomorphism. Out of this {\em 1228} are orientable and {\em 5075} are non orientable.
\end{theo}

\begin{theo}\label{thm2} For each orientable genus $g \geq 0$ there exists a centrally symmetric quadrangulated map with {\em 18g+26} vertices.
\end{theo}

\begin{theo}\label{thm3} For each orientable genus $g \geq 0$ there exists a centrally symmetric map on {\em 10g+20} vertices all of whose faces are pentagons.
\end{theo}

\begin{theo}\label{thm4} There exists a centrally symmetric orientable map of genus $g = 2k-1 (\geq 1)$ on $24 +  [\frac{2k-1}{2}] 12$ vertices all of whose faces are hexagons.
\end{theo}

\begin{theo}\label{thm5}
There are exactly {\em 68} centrally symmetric triangulated 3\linebreak-manifolds on {\em 12} vertices up to isomorphism.
\end{theo}

\section{Enumeration results for triangulated surfaces}

In this section we present an enumeration of CST surfaces with $\mathbb{Z}_{2}$ action. This is different from the vertex transitive enumeration done Lutz in \cite{lutz0}. We have modified the program MANIFOLD$_{-}$VT\cite{lutz1} of Lutz. In this program Lutz has used cyclic and dihedral group of order {\em 2m} and {\em 4m} respectively and generated CS vertex transitive triangulated surfaces whose automorphism group are $Z_{2m} = <(123\dots 2m)>$ and $D_{2m} = <(123\dots 2m), (1,~2m)(2,~2m-1)\dots(m,~m+1)>$. We have replaced the groups by $\mathbb{Z}_{2}$ and relaxed the criteria of vertex transitivity. We have used group action $\mathbb{Z}_{2}= < I >$ on the set $\{1, 2, \dots, 2m\}$ where $I = (1, 2m)(2, 2m-1)\dots(m, m+1)$ denote the generator. It generates all possible  1- and 2-orbits, that is,  1 and  2 dimensional orbits. We denote by $F^{I}$ the image of the face {\em F} under the group action $\mathbb{Z}_{2}$. We neglect those  2-orbits containing {\em F} and $F^{I}$ for which $F \bigcap F^{I} \not= \emptyset$. And we ignore those  1-orbits for which $e = e^{I}$ . The remaining  orbits are called admissible orbits. Therefore, for fixed $n = 2m$, we obtained all admissible  1- and  2-orbits under the group action $\mathbb{Z}_{2}$. In the process we check link of $ m$ vertices namely $1, 2, \dots, m$ which are use to define {\em I}. We also compute reduced homology groups to check orientability of the objects using \cite{hom}. Hence we get all possible non isomorphic CST surfaces. As a result for $m = 3, 4$ and 5 we have listed the objects in Table \ref{table1}. For $m = 3$ the object $6_{tight}$ obtained in Table \ref{table1} is isomorphic to Lutz's object \cite{lutz0}. For $ m = 4$ we get  4 objects out of which the list object $8_{tight}$ in Table \ref{table1} is isomorphic to that of Lutz's object \cite{lutz0}. For $ m = 6$, we give the number of non isomorphic objects for different genus in Table \ref{table2}. In this case and for $\chi = -8$ we give the list of all the objects in Table \ref{table3}.

We give a technique to construct a CST surface from already known CST surfaces. We use this idea in the following sections. Let $M_{1}$, $M_{2}$ be two CST surfaces and the involution $I_{1} := (a^{1}_{1}, a^{1}_{2m_{1}})$ $(a^{1}_{2}, a^{1}_{2m_{1}-1})$ \dots $(a^{1}_{m_{1}}, a^{1}_{m_{1}+1})$ act on $M_{1}$ and the involution $I_{2} := (a^{2}_{1}, a^{2}_{2m_{2}})$ $(a^{2}_{2}, a^{2}_{2m_{2}-1})\dots$ $(a^{2}_{m_{2}}, a^{2}_{m_{2}+1})$ act on $M_{2}$. We consider two faces $F_{M_{t}, 1}$ and $F_{M_{t}, 2}$ of $M_{t}$ for $t \in \{1, 2\}$ with the following properties $\colon$ $F_{M_{t},1}^{I_{t}} = F_{M_{t},2}$ and there is no edge between the vertices of $F_{M_{t},1}$ with the vertices of $F_{M_{t},2}$. We give list of maps with the following properties in {\em Section 3, 4, 5}. Let $F_{M_{1}, 1} := a^{1}_{i}a^{1}_{j}a^{1}_{k}$. Then, by the property $F_{M_{1},1}^{I_{1}} = F_{M_{1},2}$, $F_{M_{1}, 2} = a^{1}_{2m_{1}-(i-1)}$ $a^{1}_{2m_{1}-(j-1)}$ $a^{1}_{2m_{1}-(k-1)}$, see  Figure 1. Similarly, let $F_{M_{2}, 1} := a^{2}_{i}a^{2}_{j}a^{2}_{k}$ then $F_{M_{2}, 2}:= a^{2}_{2m_{2}-(i-1)}$ $a^{2}_{2m_{2}-(j-1)}$ $a^{2}_{2m_{2}-(k-1)}$, see  Figure 1. We remove interior of $F_{M_{t},1}$ and $F_{M_{t},2}$ and obtain cycles $\partial F_{M_{t},1}, \partial F_{M_{t},2}$ for $t \in \{1, 2\}$, see  Figure 1. We identify $\partial F_{M_{1},t}$ with $\partial F_{M_{2},t}$ for $t \in \{1, 2\}$ by the map $a^{1}_{s}\mapsto a^{2}_{s}$ and $a^{1}_{2m_{1}-(s-1)} \mapsto a^{2}_{2m_{2}-(s-1)}$ for $s \in \{i, j, k\}$.
%More precisely, for $t \in \{1, 2\}$, we define a map $f : \partial F_{M_{1},t} \rightarrow \partial F_{M_{2},t}$ by $f(a^{1}_{s}) = a^{2}_{s},$ $f(a^{1}_{2m_{1}-(s-1)}) = a^{2}_{2m_{2}-(s-1)},$ $f(a^{1}_{s_{1}}a^{1}_{s_{2}})=a^{2}_{s_{1}}a^{2}_{s_{2}}$ and $f(a^{1}_{2m_{1}-(s_{1}-1)}a^{1}_{2m_{1}-(s_{2}-1)})=a^{2}_{2m_{2}-(s_{1}-1)}a^{2}_{2m_{1}-(s_{2}-1)}$ for $s, s_{1}, s_{2} \in \{i, j, k\}$  and $s_{1} \not= s_{2}$, where $a^{t}_{s_{1}}a^{t}_{s_{2}}$ represent an edge.
We identify $a^{1}_{s_{1}}a^{1}_{s_{2}}$ with $f(a^{1}_{s_{1}}a^{1}_{s_{2}})$ and $a^{1}_{2m_{1}-(s_{1}-1)}a^{1}_{2m_{1}-(s_{2}-1)}$ with $f(a^{1}_{2m_{1}-(s_{1}-1)}a^{1}_{2m_{1}-(s_{2}-1)})$ for $s_{1}, s_{2} \in \{i, j, k\}$. In this process, the vertices and edges of $\partial F_{M_{1},t}$ identify with the vertices and edges of $\partial F_{M_{2},t}$ respectively. Let $\chi(M_{t})$ be Euler characteristic  of $M_{t}$ for $t \in \{1,2\}$. The Euler characteristic  of the resulting object, say {\em M}, is $\chi(M_{1}) + \chi(M_{2}) - 4$. Let {\em v} be a vertex of {\em M}. If the faces containing {\em v} are of the form $vv_{1,1}v_{1,2}, vv_{2,1}v_{2,2}, \dots, vv_{m,1}v_{m,2}$ such that $v_{m,2}= v_{1,1}$, $v_{1,2}= v_{2,1}$, \dots, $v_{m-1,2}= v_{m,1}$ for $m\geq3$ then the link {\em lk(v)} is the cycle $C(v_{1,1}, v_{1,2}, v_{2,1}, v_{2,2}, \dots, v_{m-1,1}, v_{m-1,2}, v_{m,1}, v_{m,2})$, see \cite{datta}. By the above properties, link of each vertex of $\partial F_{M_{1},1} (= \partial F_{M_{2},1}$) and $\partial F_{M_{1},2} (= \partial F_{M_{2},2}$) in {\em M} is homeomorphic to 1-sphere. Therefore {\em M} is manifold. Next, we show {\em M} is CS. We define the involution $I :=  I_{1}I_{2}-\{(a^{1}_{i}, a^{1}_{2m_{1}-(i-1)}), (a^{1}_{j}, a^{1}_{2m_{1}-(j-1)}), (a^{1}_{k}, a^{1}_{2m_{1}-(k-1)})\}$. For an example, if  $I_{1} = (1, 2) (3, 4) (5, 6)$ and $I_{2} =$(7, 8) (9, 10) (11, 12) then {\em $I=I_{1}I_{2}-\{$(3, 4) (5, 6)$\}$} gives $I = $(1, 2) (7, 8) (9, 10) (11, 12). We claim that {\em M} is CS under the involution {\em I}. Let {\em F} be a face of {\em M} such that $F \bigcap F^{I} \not= \emptyset$. Then the face {\em F} or its subset, say {\em E}, which is a face of {\em M} belongs to either one of $M_{i}$ for $ i = 1, 2$ or $M_{1}\bigcap M_{2}$. If {\em F} belongs to one of $M_{i}$ for $ i = 1 or 2$ then by the definition of involution $F \bigcap F^{I_{i}} \not=\emptyset$ and which is a contradiction. If {\em F} belongs to $M_{1}\bigcap M_{2}$ then both $I_{1}$ and $I_{2}$ fix the face {\em F}. This gives a contradiction as $M_{1}$ and $M_{2}$ are CST surfaces. We use the same above argument for {\em E}. Therefore {\em M} is CST surface.

%\newpage
%\vspace{ mm}
\begin{center}
\begin{picture}(0,0)(-5,20)
\setlength{\unitlength}{10mm}
\drawpolygon (-7,0)(-3,0)(-5,2)
\drawpolygon (.5,0)(4.5,0)(2.75,2)
%\put(2,1){$\Longrightarrow$}
%\drawpolygon (3,0)(7,0)(5,2)
\put(-7.4,-.2){\tiny $a^{1}_{i}$}
\put(-3,-.2){\tiny $a^{1}_{j}$}
\put(-4.9,2){\tiny $a^{1}_{k}$}
\put(.1,-.3){\tiny $a^{1}_{2m_{1}-(i-1)}$}
\put(4.5,-.3){\tiny $a^{1}_{2m_{1}-(j-1)}$}
\put(2.8,2.1){\tiny $a^{1}_{2m_{1}-(k-1)}$}
%\put(2.5,0){$v_{1}$}
%\put(7.25,0){$v_{2}$}
%\put(4.6,2){$v_{3}$}

\put(-5.3,.7){\tiny $F_{M_{1}, 1}$}
\put(2.25,.7){\tiny $F_{M_{1}, 2}$}

\put(-6,-1){\tiny $\partial F_{M_{1}, 1}(a^{1}_{i}, a^{1}_{j}, a^{1}_{k})$}
\put(0,-1){\tiny $\partial F_{M_{1}, 2}(a^{1}_{2m_{1}-(i-1)}, a^{1}_{2m_{1}-(j-1)}, a^{1}_{2m_{1}-(k-1)})$}
%\put(4.5,-.75){$C_{2}$}
%\put(-2, -2){$Figure-1$}

\drawpolygon (-7,-5)(-3,-5)(-5,-3)
\drawpolygon (.5,-5)(4.5,-5)(2.75,-3)
%\put(2,1){$\Longrightarrow$}
%\drawpolygon (3,0)(7,0)(5,2)
\put(-7.3,-5.2){\tiny $a^{2}_{i}$}
\put(-2.9,-5.2){\tiny $a^{2}_{j}$}
\put(-4.9,-3){\tiny $a^{2}_{k}$}
\put(.1,-5.3){\tiny $a^{2}_{2m_{2}-(i-1)}$}
\put(4.5,-5.3){\tiny $a^{2}_{2m_{2}-(j-1)}$}
\put(2.8,-3){\tiny $a^{2}_{2m_{2}-(k-1)}$}

\put(-5.3,-4.5){\tiny $F_{M_{2}, 1}$}
\put(2.25,-4.5){\tiny $F_{M_{2}, 2}$}

\put(-6,-6){\tiny $\partial F_{M_{2}, 1}(a^{2}_{i}, a^{2}_{j}, a^{2}_{k})$}
\put(0,-6){\tiny $\partial F_{M_{2}, 2}(a^{2}_{2m_{2}-(i-1)}, a^{2}_{2m_{2}-(j-1)}, a^{2}_{2m_{2}-(k-1)})$}

\put(-2, -7.5){\tiny { Figure 1}}

\end{picture}
\end{center}
\vspace{-1in}
%\hrule

\vspace{95 mm}

\vspace{1in}

\begin{example}\label{example1}
{\small {\bf Construction of CST orientable surface with $\chi = -4$ from already known CST torus }}
\end{example}

Let $M_{1} := \{$ $123,$ $124,$ $134,$ $235,$ $245,$ $346,$ $357,$ $369,$ $378,$ $389,$ $450^{'},$ $468,$ $478,$ $470^{'},$ $569,$ $560^{'},$ $579,$ $680^{'},$ $790^{'},$ $891^{'},$ $80^{'}1^{'},$ $90^{'}2^{'},$ $91^{'}2^{'},$ $0^{'}1^{'}2^{'}\}$ be a triangulation of torus. The map $M_{1}$ is CS under the involution $I_{1} = (1, 2^{'})$ $(2, 1^{'})$ $(3, 0^{'})$ $(4, 9)$ $(5, 8)$ $(6, 7)$. We consider an isomorphic copy of $M_{1}$ and denote it by $M_{2}$. Let $M_{2} := \{$ $3^{'}4^{'}5^{'},$ $3^{'}4^{'}6^{'},$ $3^{'}5^{'}6^{'},$ $4^{'}5^{'}7^{'},$ $4^{'}6^{'}7^{'},$ $5^{'}6^{'}8^{'},$ $5^{'}7^{'}9^{'},$ $5^{'}8^{'}1^{''},$ $5^{'}9^{'}0^{''},$ $5^{'}0^{''}1^{''},$ $6^{'}7^{'}2^{''},$ $6^{'}8^{'}0^{''},$ $6^{'}9^{'}0^{''},$ $6^{'}9^{'} 2^{''},$ $7^{'}8^{'}1^{''},$ $7^{'}8^{'}2^{''},$ $7^{'}9^{'}1^{''},$ $8^{'}0^{''}2^{''},$ $9^{'}1^{''}2^{''},$ $0^{''}1^{''}3^{''},$ $ 0^{''}2^{''}3^{''},$ $1^{''}2^{''}4^{''},$ $1^{''}3^{''}4^{''},$ $2^{''}3^{''}4^{''}\}$. The map $M_{2}$ is CS under the involution $I_{2} = (3^{'}, 4^{''})$ $(4^{'}, 3^{''})$ $(5^{'}, 2^{''})$ $(6^{'}, 1^{''})$ $(7^{'}, 0^{''})$ $(8^{'}, 9^{'})$. We remove interior of two faces namely $123,$ $0^{'}1^{'}2^{'}$ from $M_{1}$ and $3^{'}4^{'}5^{'},$ $2^{''}3^{''}4^{''}$ from $M_{2}$. We identify the cycles $C(1, 2, 3)$ with $C(3^{'}, 4^{'}, 5^{'})$ and $C(0^{'}, 1^{'}, 2^{'})$ with $C(2^{''}, 3^{''}, 4^{''})$ by the map $1 \mapsto 3^{'},$ $2 \mapsto 4^{'},$ $3 \mapsto 5^{'},$ $0^{'} \mapsto 2^{''},$ $1^{'} \mapsto 3^{''},$ $2^{'} \mapsto 4^{''}$. So, we get the map $M : = \{$ $ 124,$ $134,$ $235,$ $245,$ $346,$ $357,$ $369,$ $378,$ $389,$ $450^{'},$ $468,$ $478,$ $470^{'},$ $569,$ $560^{'},$ $579,$ $680^{'},$ $790^{'},$ $891^{'},$ $80^{'}1^{'},$ $90^{'}2^{'},$ $91^{'}2^{'},$ $126^{'},$ $136^{'},$ $237^{'},$ $26^{'}7^{'},$ $36^{'}8^{'},$ $37^{'}9^{'},$ $38^{'}1^{''},$ $39^{'}0^{''},$ $30^{''}1^{''},$ $ 6^{'}7^{'}0^{'},$ $6^{'}8^{'}0^{''},$ $6^{'}9^{'}0^{''},$ $6^{'}9^{'}0^{'},$ $7^{'}8^{'}1^{''},$ $7^{'}8^{'}0^{'},$ $7^{'}9^{'}1^{''},$ $8^{'}0^{''}0^{'},$ $9^{'}1^{''}0^{'},$ $0^{''}1^{''}1^{'},$ $0^{''}0^{'}1^{'},$ $ 1^{''}0^{'}2^{'},$ $1^{''}1^{'}2^{'}\}$. The Euler characteristic  $\chi(M)$ is $-4$. The map {\em M} is CS under the involution $I = (1, 2^{'})$ $(2, 1^{'})$ $(3, 0^{'})$ $(4, 9)$ $(5, 8)$ $(6, 7)$ $(6^{'}, 1^{''})$ $(7^{'}, 0^{''})$ $(8^{'}, 9^{'})$.

\bigskip

\noindent{\sc Proof of Theorem}\ref{thm1} The proof follows from arguments in  {\em Section 2}. Table \ref{table1} gives the list of centrally symmetric triangulated surfaces for { n $\leq$ 10} vertices. Table \ref{table2} gives number of different objects on { 12} vertices. The total number of objects is { 6303}. It is clear from the tables by looking at homology groups that { 1228} are orientable and { 5075} are non orientable.
\hfill$\Box$

\begin{remark}\label{remark1}
As one can see from Table \ref{table3} the number of CS objects on { 12} vertices is huge and it is practically not possible to give all the list here. However we give the list of objects with $\chi = -8$. Other objects are available with authors and may be supplied on demand. For the sake of clarity, in Table \ref{table3}, we have used $a$ and $b$ to denote 10 and 11 respectively. In each table we have used {\em ijk} to represent a 2-orbit in the list of orbits of order 2, where $i, j, k \in \{1, 2, \dots, 8, 9, a, b\}$.
\end{remark}

\normalsize
%\newpage
\begin{table}
\tiny
%\begin{table}[ht]
\caption{List of the CST Surfaces}
\centering % used for centering table
\begin{tabular}{c c c c c} % centered columns (4 columns)
\hline\hline %inserts double horizontal lines
$n$ & $2$-Manifold & $f$-vector & List of Orbits  \\ [0.3ex] % inserts table
%heading
\hline % inserts single horizontal line
$6_{tight}$ & $S^{2}$ & ( 6, 12, 8 )  & 123, 124, 135, 145\\

8 & $S^{2}$ & (8, 18, 12) &  123, 124, 134, 235, 246, 256\\

8 & $S^{2}$ & (8, 18, 12) &  123, 124, 135, 146, 156, 234\\

8 & $S^{2}$ & (8, 18, 12) & 123, 124, 135, 146, 157, 167\\

8 & $S^{2}$ & (8, 18, 12) & 123, 124, 137, 147, 235, 246\\

$8_{tight}$ & $S^{1}\times S^{1}$ & (8, 24, 16) & 123, 124, 135, 147, 156, 167, 246, 256\\

10 & $S^{2}$ & (10, 24, 16)  & 123, 124, 134, 235, 245, 346, 357, 367\\

10 & $S^{2}$ & (10, 24, 16) &  123, 124, 134, 235, 246, 257, 268, 278\\

10 & $S^{2}$ & (10, 24, 16) & 123, 124, 134, 235, 248, 257, 278, 357\\

10 & $S^{2}$ & (10, 24, 16) & 123, 124, 135, 145, 234, 346, 357, 367\\

10 & $S^{2}$ & (10, 24, 16) & 123, 124, 135, 145, 236, 246, 357, 367\\

10 & $S^{2}$ & (10, 24, 16) & 123, 124, 135, 146, 157, 167, 234, 345\\

10 & $S^{2}$ & (10, 24, 16) & 123, 124, 135, 146, 157, 168, 178, 234\\

10 & $S^{2}$ & (10, 24, 16) & 123, 124, 135, 146, 157, 168, 179, 189\\

10 & $S^{2}$ & (10, 24, 16) & 123, 124, 135, 146, 159, 169, 235, 246\\

10 & $S^{2}$ & (10, 24, 16) & 123, 124, 135, 148, 158, 234, 346, 357\\

10 & $S^{2}$ & (10, 24, 16) & 123, 124, 135, 148, 158, 236, 246, 357\\

10 & $S^{2}$ & (10, 24, 16) & 123, 124, 135, 148, 159, 189, 246, 357\\

10 & $S^{2}$ & (10, 24, 16) & 123, 124, 135, 148, 159, 189, 248, 268\\

10 & $S^{2}$ & (10, 24, 16) & 123, 124, 137, 148, 178, 235, 245, 345\\

10 & $S^{2}$ & (10, 24, 16) & 123, 124, 139, 149, 235, 246, 257, 268\\

10 & $S^{2}$ & (10, 24, 16) & 123, 124, 139, 149, 235, 248, 257, 357\\

10 & $S^{1}\times S^{1}$ & (10, 30, 20) & 123, 124, 134, 235, 246, 258, 267, 278, 357, 367\\

10 & $S^{1}\times S^{1}$ & (10, 30, 20) & 123, 124, 135, 145, 236, 248, 257, 258, 267, 357 \\

10 & $S^{1}\times S^{1}$ & (10, 30, 20) & 123, 124, 135, 146, 157, 167, 236, 245, 258, 268\\

10 & $S^{1}\times S^{1}$ & (10, 30, 20) & 123, 124, 135, 146, 157, 167, 236, 248, 268, 367\\

10 & $S^{1}\times S^{1}$ & (10, 30, 20) & 123, 124, 135, 146, 157, 168, 178, 237, 248, 278\\

10 & $S^{1}\times S^{1}$ & (10, 30, 20) & 123, 124, 135, 146, 158, 167, 178, 234, 357, 367\\

10 & $S^{1}\times S^{1}$ & (10, 30, 20) & 123, 124, 135, 146, 158, 167, 179, 189, 357, 367\\

10 & $S^{1}\times S^{1}$ & (10, 30, 20) & 123, 124, 135, 146, 158, 169, 189, 245, 357, 367\\

10 & $S^{1}\times S^{1}$ & (10, 30, 20) & 123, 124, 135, 146, 158, 169, 189, 246, 258, 268\\

10 & $S^{1}\times S^{1}$ & (10, 30, 20) & 123, 124, 135, 146, 158, 169, 189, 248, 258, 357\\

10 & $S^{1}\times S^{1}$ & (10, 30, 20) & 123, 124, 135, 146, 159, 169, 236, 245, 357, 367\\

10 & $S^{1}\times S^{1}$ & (10, 30, 20) & 123, 124, 135, 146, 159, 169, 236, 248, 258, 357\\

10 & $S^{1}\times S^{1}$ & (10, 30, 20) & 123, 124, 135, 148, 158, 236, 245, 257, 267, 357\\

10 & $S^{1}\times S^{1}$ & (10, 30, 20) & 123, 124, 135, 148, 158, 237, 246, 257, 258, 268\\

10 & $S^{1}\times S^{1}$ & (10, 30, 20) & 123, 124, 135, 148, 158, 237, 246, 268, 278, 346\\

10 & $S^{1}\times S^{1}$ & (10, 30, 20) & 123, 124, 135, 148, 159, 167, 168, 179, 236, 367\\

10 & $S^{1}\times S^{1}$ & (10, 30, 20) & 123, 124, 135, 148, 159, 189, 245, 257, 267, 357\\

10 & $S^{1}\times S^{1}$ & (10, 30, 20) & 123, 124, 135, 149, 158, 167, 168, 179, 236, 267\\

10 & $S^{1}\times S^{1}$ & (10, 30, 20) & 123, 124, 135, 149, 158, 168, 169, 236, 245, 267\\

10 & $S^{1}\times S^{1}$ & (10, 30, 20) & 123, 124, 135, 149, 158, 168, 169, 237, 245, 367\\

10 & $S^{1}\times S^{1}$ & (10, 30, 20) & 123, 124, 135, 149, 158, 189, 245, 257, 357, 367\\

10 & $S^{1}\times S^{1}$ & (10, 30, 20) & 123, 124, 135, 149, 158, 189, 246, 257, 258, 268\\

10 & $S^{1}\times S^{1}$ & (10, 30, 20) & 123, 124, 135, 149, 158, 189, 246, 268, 278, 346\\

10 & $S^{1}\times S^{1}$ & (10, 30, 20) & 123, 124, 137, 148, 178, 235, 246, 257, 267, 345\\

10 & $S^{1}\times S^{1}$ & (10, 30, 20) & 123, 124, 137, 148, 178, 235, 246, 257, 268, 278\\

10 & $S^{1}\times S^{1}$ & (10, 30, 20) & 123, 124, 137, 149, 178, 189, 245, 257, 346, 367\\

10 & $S^{1}\times S^{1}$ & (10, 30, 20) & 123, 124, 137, 149, 178, 189, 245, 258, 278, 367\\

10 & $S^{1}\times S^{1}$ & (10, 30, 20) & 123, 124, 139, 149, 235, 246, 258, 267, 357, 367\\

10 & $S^{1}\times S^{1}$ & (10, 30, 20) & 123, 124, 139, 149, 235, 248, 257, 345, 346, 367\\

10 & $Klein~bottle$ & (10, 30, 20) & 123, 124, 135, 146, 157, 169, 178, 189, 245, 345\\

10 & $Klein~bottle$ & (10, 30, 20) & 123, 124, 135, 146, 158, 167, 179, 189, 345, 346\\

10 & $Klein~bottle$ & (10, 30, 20) & 123, 124, 135, 146, 158, 169, 189, 245, 345, 346\\

10 & $Klein~bottle$ & (10, 30, 20) & 123, 124, 135, 146, 159, 169, 236, 245, 345, 346\\

10 & $Klein~bottle$ & (10, 30, 20) & 123, 124, 135, 148, 158, 237, 245, 257, 345, 346\\

10 & $Klein~bottle$ & (10, 30, 20) & 123, 124, 135, 148, 158, 237, 245, 258, 278, 345\\

10 & $Klein~bottle$ & (10, 30, 20) & 123, 124, 135, 149, 158, 189, 245, 257, 345, 346\\

10 & $Klein~bottle$ & (10, 30, 20) & 123, 124, 135, 149, 158, 189, 245, 258, 267, 268\\

10 & $Klein~bottle$ & (10, 30, 20) & 123, 124, 135, 149, 158, 189, 245, 258, 278, 345\\

10 & $Klein~bottle$ & (10, 30, 20) & 123, 124, 137, 149, 178, 189, 245, 257, 345, 357\\

10 & $Klein~bottle$ & (10, 30, 20) & 123, 124, 139, 149, 235, 246, 258, 267, 345, 346\\
[1ex] % [1ex] adds vertical space
\hline %inserts single line
\end{tabular}
\label{table1} % is used to refer this table in the text
\end{table}

\begin{table}
\tiny
\caption{Number of CST maps on 12 vertices} % title of Table
\centering % used for centering table
\begin{tabular}{c c} % centered columns (4 columns)
\hline\hline %inserts double horizontal lines
Homology Groups (H$_{0}$, H$_{1}$, H$_{2}$) &  Number of different objects\\ [0.1ex]
\hline % inserts single horizontal line
(1, 0, 1) & 81 \\
(1, 2, 1) & 499 \\
(1, 1+Z$_{2}$, 0) & 232 \\
(1, 4, 1) & 178 \\
(1, 3+Z$_{2}$, 0 ) & 1180 \\
(1, 6, 1) & 154 \\
(1, 5+Z$_{2}$, 0) & 2707 \\
(1, 8, 1) & 258 \\
(1, 7+Z$_{2}$, 0) & 918 \\
(1, 10, 1) & 7 \\
(1, 9+Z$_{2}$, 0) & 27\\
[1ex] % [1ex] adds vertical space
\hline %inserts single line
\end{tabular}
\label{table:nonlin} % is used to refer this table in the text
\label{table2}
\end{table}
%\normalsize

\begin{table}
\tiny
\caption{List of CST maps on 12 vertices with $\chi=-8$} % title of Table
\centering % used for centering table
\begin{tabular}{c c} % centered columns (4 columns)
\hline\hline %inserts double horizontal lines
Homology Groups (H$_{0}$, H$_{1}$, H$_{2}$) & Orbits of Object\\ [0.1ex]
\hline % inserts single horizontal line
(1, 10, 1) & 123, 124, 135, 146, 157, 16a, 17b, 189, 18b, 19a, 236, 247,259, 278, 28a, 29a, 369, 378, 389, 468\\

(1, 10, 1) & 123, 124, 135, 146 157, 16a, 17b, 189, 18b, 19a, 236, 248, 257, 279, 28a, 29a, 369, 378, 389, 478\\

(1, 10, 1) & 123, 124, 135, 146, 157, 189, 238, 257, 269, 289, 346, 359, 378, 478, 16a, 17b, 18b, 19a, 24a, 27a\\

(1, 10, 1) & 123, 124, 135, 146, 159, 178, 236, 247, 268, 289, 357, 379, 389, 478, 16b, 17a, 18b, 19a, 25a, 29a\\

(1, 10, 1) & 123, 124, 135, 146, 159, 178, 239, 247, 268, 289, 345, 368, 379, 478, 16b, 17a, 18b, 19a, 25a, 26a\\

(1, 10, 1) &  123, 124, 135, 146, 159, 189, 238, 257, 259, 268, 347, 369, 378, 456, 16b, 17a, 17b, 18a, 24a, 29a\\

(1, 10, 1) &  123, 124, 135, 146, 159, 189, 239, 248, 256, 278, 347, 368, 379, 457, 16b, 17a, 17b, 18a, 25a, 29a\\

(1, 9+Z$_{2}$, 0) & 123, 124, 135, 146, 178, 189, 236, 248, 257, 259, 369, 378, 379, 457, 15a, 16b, 17b, 19a, 28a, 29a\\

(1, 9+Z$_{2}$, 0) & 123, 124, 135, 146, 178, 189, 238, 257, 259, 289, 346, 357, 369, 478, 15a, 16b, 17b, 19a, 24a, 26a\\

(1, 9+Z$_{2}$, 0) & 123, 124, 135, 148, 169, 179, 237, 257, 268, 289, 348, 356, 369, 456, 15a, 16a, 17b, 18b, 24a, 29a\\

(1, 9+Z$_{2}$, 0) &  123, 124, 135, 148, 169, 189, 239, 256, 278, 289, 346, 347, 357, 468, 15a, 16b, 17a, 17b, 24a, 25a\\

(1, 9+Z$_{2}$, 0) & 123, 124, 135, 168, 169, 178, 237, 245, 279, 289, 346, 359, 368, 478, 14a, 15b, 17b, 19a, 25a, 26a\\

(1, 9+Z$_{2}$, 0) & 123, 124, 135, 168, 169, 178, 238, 245, 259, 279, 345, 368, 379, 478, 14a, 15b, 17b, 19a, 26a, 27a\\

(1, 9+Z$_{2}$, 0) & 123, 124, 135, 168, 179, 189, 236, 248, 256, 279, 347, 357, 389, 478, 14a, 15b, 16b, 17a, 25a, 29a\\

(1, 9+Z$_{2}$, 0) & 123, 124, 135, 148, 169, 179, 237, 257, 268, 289, 346, 356, 389, 478, 15a, 16a, 17b, 18b, 24a, 29a\\

(1, 9+Z$_{2}$, 0) & 123, 124, 135, 148, 169, 179, 239, 245, 268, 278, 345, 368, 379, 478, 15a, 16a, 17b, 18b, 27a, 29a\\

(1, 9+Z$_{2}$, 0) & 123, 124, 135, 146, 159, 178, 238, 257, 259, 289, 347, 368, 369, 456, 16b, 17b, 18a, 19a, 24a, 26a\\

(1, 9+Z$_{2}$, 0) & 123, 124, 135, 146, 159, 178, 236, 248, 257, 289, 369, 378, 379, 457, 16b, 17b, 18a, 19a, 25a, 29a\\

(1, 9+Z$_{2}$, 0) & 123, 124, 135, 146, 159, 178, 239, 245, 268, 289, 347, 369, 378, 468, 16b, 17b, 18a, 19a, 25a, 27a\\

(1, 9+Z$_{2}$, 0) & 123, 124, 135, 146, 159, 178, 239, 245, 278, 289, 347, 368, 369, 468, 16b, 17b, 18a, 19a, 25a, 26a\\

(1, 9+Z$_{2}$, 0) & 123, 124, 135, 146, 157, 169, 178, 247, 259, 268, 289, 348, 356, 379, 389, 18b, 19a, 26a, 27a, 1ab\\

(1, 9+Z$_{2}$, 0) & 123, 124, 135, 146, 157, 169, 189, 247, 256, 268, 289, 348, 359, 378, 379, 17a, 18b, 27a, 29a, 1ab\\

(1, 9+Z$_{2}$, 0) & 123, 124, 135, 146, 157, 169, 189, 247, 259, 268, 278, 346, 359, 378, 389, 17a, 18b, 26a, 29a, 1ab\\

(1, 9+Z$_{2}$, 0) & 123, 124, 135, 146, 157, 169, 189, 259, 268, 278, 279, 345, 348, 369, 378, 17a, 18b, 24a, 26a, 1ab\\

(1, 9+Z$_{2}$, 0) & 123, 124, 135, 146, 157, 169, 237, 256, 268, 289, 348, 379, 389, 478, 17a, 18a, 18b, 19b, 27a, 29a\\

(1, 9+Z$_{2}$, 0) & 123, 124, 135, 146, 157, 189, 238, 256, 259, 278, 347, 369, 389, 468, 16a, 17b, 18a, 19b, 27a, 29a\\

(1, 9+Z$_{2}$, 0) & 123, 124, 135, 146, 157, 189, 238, 259, 278, 279, 348, 356, 369, 457, 16a, 17b, 18b, 19a, 24a, 26a\\

(1, 9+Z$_{2}$, 0) & 123, 124, 135, 146, 157, 189, 238, 256, 257, 289, 348, 369, 379, 478, 16b, 17a, 18a, 19b, 26a, 29a\\

(1, 9+Z$_{2}$, 0) & 123, 124, 135, 146, 159, 178, 179, 257, 269, 278, 289, 345, 348, 368, 369, 16a, 18b, 24a, 26a, 1ab\\

(1, 9+Z$_{2}$, 0) & 123, 124, 135, 146, 159, 178, 239, 257, 268, 269, 347, 348, 389, 456, 16a, 17a, 18b, 19b, 27a, 28a\\

(1, 9+Z$_{2}$, 0) & 123, 124, 135, 146, 159, 178, 239, 268, 278, 279, 345, 348, 369, 457, 16a, 17a, 18b, 19b, 25a, 26a\\

(1, 9+Z$_{2}$, 0) & 123, 124, 135, 146, 159, 189, 238, 257, 268, 279, 345, 369, 378, 478, 16a, 17a, 17b, 18b, 24a, 29a\\

(1, 9+Z$_{2}$, 0) & 123, 124, 135, 146, 159, 189, 238, 257, 269, 278, 346, 357, 389, 478, 16a, 17a, 17b, 18b, 24a, 29a\\

(1, 9+Z$_{2}$, 0) & 123, 124, 135, 146, 159, 178, 189, 248, 257, 268, 269, 345, 369, 378, 379, 16b, 17a, 25a, 29a, 1ab\\
[1ex] % [1ex] adds vertical space
\hline %inserts single line
\end{tabular}
\label{table:nonlin} % is used to refer this table in the text
\label{table3}
\end{table}
%\normalsize

\normalsize
%\newpage
\begin{table}
\tiny
%\begin{table}[ht]
\caption{List of the CST { 3}-manifolds on $12$ vertices with reduced homology groups ($H_{0}, H_{1}, H_{2}, H_{3}$) = (1, 1, 1, 1).}
\centering % used for centering table
\begin{tabular}{c c} % centered columns (4 columns)
\hline\hline %inserts double horizontal lines
Sl. No. & List of Orbits  \\ [0.3ex] % inserts table
%heading
\hline % inserts single horizontal line
$M_{1}$ &   1234, 1235, 1245, 1348, 1356, 1368, 1456, 1689, 2347, 2457, 2579, 3456, 3468, 146b, 148b, 169a, 189a, 16ab, 18ab\\

$M_{2}$ &   1234, 1235, 1245, 1348, 1356, 1368, 1456, 1689, 2347, 2457, 3456, 146b, 148b, 169a, 189a, 257a, 259a, 279a, 16ab, 18ab\\

$M_{3}$ &   1234, 1235, 1245, 1348, 1356, 1368, 1456, 1689, 2357, 2579,  3456, 3457, 3468, 146b, 148b, 169a, 189b, 16ab, 19ab\\

$M_{4}$ &   1234, 1235, 1245, 1348, 1356, 1368, 1456, 1689, 2357, 3456, 3457, 146b, 148b, 169a, 189b, 257a, 259a, 279a, 16ab, 19ab\\

$M_{5}$ &  1234, 1235, 1245, 1348, 1356, 1368, 1456, 1689, 2457, 2579, 3456, 3468, 146b, 148b, 169b, 189a, 18ab, 19ab\\

$M_{6}$ &  1234, 1235, 1245, 1348, 1356, 1368, 1456, 1689, 2457, 3456, 146b, 148b, 169b, 189a, 257a, 259a, 279a, 18ab, 19ab\\

$M_{7}$ &  1234, 1235, 1245, 1348, 1356, 1368, 1456, 1689, 2347, 2357, 2579, 3456, 3457, 3468, 146b, 148b, 169b, 189b\\

$M_{8}$ &  1234, 1235, 1245, 1348, 1356, 1368, 1456, 1689, 2347, 146b, 148b, 169b, 189b, 257a, 259a, 279a\\

$M_{9}$ &   1234, 1235, 1245, 1348, 1356, 1368, 1456, 2347, 2357, 2579, 3456, 3468, 146b, 148b, 168a, 169a, 169b, 189a, 189b\\

$M_{10}$ &  1234, 1235, 1245, 1348, 1356, 1368, 1456, 2347, 2357, 3456, 146b, 148b, 168a, 169a, 169b, 189a, 189b, 257a, 259a, 279a\\

$M_{11}$ &  1234, 1235, 1245, 1348, 1356, 1368, 1456, 2457, 2579, 3456, 3457, 3468, 146b, 148b, 168a, 169a, 169b, 18ab, 19ab\\

$M_{12}$ &  1234, 1235, 1245, 1348, 1356, 1368, 1456, 2457, 3456, 3457, 146b, 148b, 168a, 169a, 169b, 257a, 259a, 279a, 18ab, 19ab \\

$M_{13}$ &   1234, 1235, 1245, 1348, 1356, 1368, 1456, 2357, 2579, 3456, 3468, 146b, 148b, 168a, 189a, 189b, 16ab, 19ab\\

$M_{14}$ &   1234, 1235, 1245, 1348, 1356, 1368, 1456, 2357, 3456, 146b, 148b, 168a, 189a, 189b, 257a, 259a, 279a, 16ab, 19ab\\

$M_{15}$ &  1234, 1235, 1245, 1348, 1356, 1368, 1456, 2347, 2457, 2579, 3456, 3457, 3468, 146b, 148b, 168a, 16ab, 18ab\\

$M_{16}$ &   1234, 1235, 1245, 1348, 1356, 1368, 1456, 2347, 2457, 3456, 3457, 146b, 148b, 168a, 257a, 259a, 279a, 16ab, 18ab\\

$M_{17}$ &   1234, 1235, 1245, 1356, 1368, 1456, 1689, 2347, 2457, 2579, 3456, 3468, 134b, 138b, 146b, 169a, 189a, 259a, 16ab, 18ab\\

$M_{18}$ &   1234, 1235, 1245, 1356, 1368, 1456, 1689, 2347, 2457, 3456, 134b, 138b, 146b, 169a, 189a, 257a, 279a, 16ab, 18ab\\

$M_{19}$ &   1234, 1235, 1245, 1356, 1368, 1456, 1689, 2357, 2579, 3456, 3457, 3468, 134b, 138b, 146b, 169a, 189b, 259a, 16ab, 19ab\\

$M_{20}$ &   1234, 1235, 1245, 1356, 1368, 1456, 1689, 2357, 3456, 3457, 134b, 138b, 146b, 169a, 189b, 257a, 279a, 16ab, 19ab\\

$M_{21}$ &   1234, 1235, 1245, 1356, 1368, 1456, 1689, 2457, 2579, 3456,  3468, 134b, 138b, 146b, 169b, 189a, 259a, 18ab, 19ab\\

$M_{22}$ &   1234, 1235, 1245, 1356, 1368, 1456, 1689, 2457, 3456, 134b, 138b, 146b, 169b, 189a, 257a, 279a, 18ab, 19ab\\

$M_{23}$ &   1234, 1235, 1245, 1356, 1368, 1456, 1689, 2347, 2357, 2579, 3456, 3457, 3468, 134b, 138b, 146b, 169b, 189b, 259a\\

$M_{24}$ &  1234, 1235, 1245, 1356, 1368, 1456, 1689, 2347, 2357, 3456, 3457, 134b, 138b, 146b, 169b, 189b, 257a, 279a\\

$M_{25}$ &  1234, 1235, 1245, 1356, 1368, 1456, 2347, 2357, 2579, 3456, 3468, 134b, 138b, 146b, 168a, 169a, 169b, 189a, 189b, 259a\\

$M_{26}$ &   1234, 1235, 1245, 1356, 1368, 1456, 2347, 2357, 3456, 134b, 138b, 146b, 168a, 169a, 169b, 189a, 189b, 257a, 279a\\

$M_{27}$ &   1234, 1235, 1245, 1356, 1368, 1456, 2457, 2579, 3456, 3457, 3468, 134b, 138b, 146b, 168a, 169a, 169b, 259a, 18ab, 19ab\\

$M_{28}$ &   1234, 1235, 1245, 1356, 1368, 1456, 2457, 3456, 3457, 134b, 138b, 146b, 168a, 169a, 169b, 257a, 279a, 18ab, 19ab\\

$M_{29}$ &   1234, 1235, 1245, 1356, 1368, 1456, 2347, 2457, 2579, 3456, 3457, 3468, 134b, 138b, 146b, 168a, 259a, 16ab, 18ab\\

$M_{30}$ &   1234, 1235, 1245, 1356, 1368, 1456, 2347, 2457, 3456, 3457, 134b, 138b, 146b, 168a, 257a, 279a, 16ab, 18ab\\

$M_{31}$ &   1234, 1235, 1246, 1256, 1348, 1356, 1368, 1689, 2347, 2456, 2457, 2579, 3456, 3468, 146b, 148b, 169a, 189a, 16ab, 18ab\\

$M_{32}$ &   1234, 1235, 1246, 1256, 1348, 1356, 1368, 1689, 2347, 2456, 2457, 3456, 146b, 148b, 169a, 189a, 257a, 259a, 279a, 16ab, 18ab\\

$M_{33}$ &   1234, 1235, 1246, 1256, 1348, 1356, 1368, 1689, 2357, 2456, 2579, 3456, 3457, 3468, 146b, 148b, 169a, 189b, 16ab, 19ab\\

$M_{34}$ &   1234, 1235, 1246, 1256, 1348, 1356, 1368, 1689, 2357, 2456, 3456, 3457, 146b, 148b, 169a, 189b, 257a, 259a, 279a, 16ab, 19ab\\

$M_{35}$ &   1234, 1235, 1246, 1256, 1348, 1356, 1368, 1689, 2347, 2357, 2456, 2579, 3456, 3457, 3468, 146b, 148b, 169b, 189b\\

$M_{36}$ &   1234, 1235, 1246, 1256, 1348, 1356, 1368, 1689, 2347, 2357, 2456, 3456, 3457, 146b, 148b, 169b, 189b, 257a, 259a, 279a\\

$M_{37}$ &   1234, 1235, 1246, 1256, 1348, 1356, 1368, 2347, 2357, 2456, 2579, 3456, 3468, 146b, 148b, 168a, 169a, 169b, 189a, 189b\\

$M_{38}$ &   1234, 1235, 1246, 1256, 1348, 1356, 1368, 2347, 2357, 2456, 3456, 146b, 148b, 168a, 169a, 169b, 189a, 189b, 257a, 259a, 279a\\

$M_{39}$ &   1234, 1235, 1246, 1256, 1348, 1356, 1368, 2456, 2457, 2579, 3456, 3457, 3468, 146b, 148b, 168a, 169a, 169b, 18ab, 19ab\\

$M_{40}$ &   1234, 1235, 1246, 1256, 1348, 1356, 1368, 2456, 2457, 3456, 3457, 146b, 148b, 168a, 169a, 169b, 257a, 259a, 279a, 18ab, 19ab\\

$M_{41}$ &   1234, 1235, 1246, 1256, 1348, 1356, 1368, 2347, 2456, 2457, 2579, 3456, 3457, 3468, 146b, 148b, 168a, 16ab, 18ab\\

$M_{42}$ &  1234, 1235, 1246, 1256, 1348, 1356, 1368, 2347, 2456, 2457, 3456, 3457, 146b, 148b, 168a, 257a, 259a, 279a, 16ab, 18ab\\

$M_{43}$ &   1234, 1235, 1246, 1256, 1348, 1356, 1468, 1689, 2347, 2456, 2457, 2579, 3456, 136b, 138b, 169a, 189a, 259a, 279a, 16ab, 18ab\\

$M_{44}$ &   1234, 1235, 1246, 1256, 1348, 1356, 1468, 1689, 2347, 2456, 2457, 3456, 3468, 136b, 138b, 169a, 189a, 257a, 16ab, 18ab\\

$M_{45}$ &   1234, 1235, 1246, 1256, 1348, 1356, 1468, 1689, 2347, 2357, 2456, 2579, 3456, 3457, 136b, 138b, 169b, 189b, 259a, 279a\\

$M_{46}$ &   1234, 1235, 1246, 1256, 1348, 1356, 1468, 1689, 2347, 2357, 2456, 3456, 3457, 3468, 136b, 138b, 169b, 189b, 257a\\

$M_{47}$ &   1234, 1235, 1246, 1256, 1348, 1356, 1468, 2347, 2357, 2456, 2579, 3456, 136b, 138b, 168a, 169a, 169b, 189a, 189b, 259a, 279a\\

$M_{48}$ &   1234, 1235, 1246, 1256, 1348, 1356, 1468, 2347, 2357, 2456, 3456, 3468, 136b, 138b, 168a, 169a, 169b, 189a, 189b, 257a\\

$M_{49}$ &   1234, 1235, 1246, 1256, 1348, 1356, 1468, 2456, 2457, 3456, 3457, 3468, 136b, 138b, 168a, 169a, 169b, 257a, 18ab, 19ab\\

$M_{50}$ &   1234, 1235, 1246, 1256, 1348, 1356, 1468, 2347, 2456, 2457, 2579, 3456, 3457, 136b, 138b, 168a, 259a, 279a, 16ab, 18ab\\

$M_{51}$ &   1234, 1235, 1246, 1256, 1348, 1356, 1468, 2347, 2456, 2457, 3456, 3457, 3468, 136b, 138b, 168a, 257a, 16ab, 18ab\\

$M_{52}$ &   1234, 1235, 1246, 1256, 1356, 1368, 1689, 2347, 2456, 2457, 3456, 134b, 138b, 146b, 169a, 189a, 257a, 279a, 16ab, 18ab\\

$M_{53}$ &   1234, 1235, 1246, 1256, 1356, 1368, 1689, 2347, 2357, 2456, 2579, 3456, 3457, 3468, 134b, 138b, 146b, 169b, 189b, 259a\\

$M_{54}$ &   1234, 1235, 1246, 1256, 1356, 1368, 1689, 2347, 2357, 2456, 3456, 3457, 134b, 138b, 146b, 169b, 189b, 257a, 279a\\

$M_{55}$ &   1234, 1235, 1246, 1256, 1356, 1368, 2456, 2457, 3456, 3457, 134b, 138b, 146b, 168a, 169a, 169b, 257a, 279a, 18ab, 19ab\\

$M_{56}$ &   1234, 1235, 1246, 1256, 1356, 1368, 2347, 2456, 2457, 2579, 3456, 3457, 3468, 134b, 138b, 146b, 168a, 259a, 16ab, 18ab\\

$M_{57}$ &   1234, 1235, 1246, 1256, 1356, 1368, 2347, 2456, 2457, 3456, 3457, 134b, 138b, 146b, 168a, 257a, 279a, 16ab, 18ab\\

$M_{58}$ &   1234, 1235, 1246, 1256, 1356, 1468, 1689, 2347, 2357, 2456, 2579, 3456, 3457, 134b, 136b, 148b, 169b, 189b, 279a\\

$M_{59}$ &   1234, 1235, 1246, 1256, 1356, 1468, 1689, 2347, 2357, 2456, 3456, 3457, 3468, 134b, 136b, 148b, 169b, 189b, 257a, 259a\\

$M_{60}$ &   1234, 1235, 1246, 1256, 1356, 1468, 2347, 2357, 2456, 2579, 3456, 134b, 136b, 148b, 168a, 169a, 169b, 189a, 189b, 279a\\

$M_{61}$ &   1234, 1235, 1246, 1256, 1356, 1468, 2347, 2456, 2457, 2579, 3456, 3457, 134b, 136b, 148b, 168a, 279a, 16ab, 18ab\\

$M_{62}$ &   1234, 1235, 1246, 1256, 1356, 1468, 2347, 2456, 2457, 3456, 3457, 3468, 134b, 136b, 148b, 168a, 257a, 259a, 16ab, 18ab\\

$M_{63}$ &   1234, 1235, 1246, 1259, 1269, 1348, 1356, 1368, 1468, 1569, 2345, 2456, 2579, 3456, 3478, 256a, 257a, 269a, 279a\\

$M_{64}$ &   1234, 1235, 1246, 1259, 1269, 1348, 1356, 1368, 1468, 1569, 2346, 2356, 2579, 3478, 256a, 257a, 269a, 279a\\

$M_{65}$ &   1234, 1235, 1246, 1259, 1269, 1356, 1569, 2345, 2456, 2579, 3456, 3468, 3478, 134b, 136b, 146b, 256a, 257a, 269a\\

$M_{66}$ &   1234, 1235, 1246, 1259, 1269, 1356, 1569, 2346, 2356, 2579, 3468, 3478, 134b, 136b, 146b, 256a, 257a, 269a\\

$M_{67}$ &   1234, 1235, 1246, 1356, 2345, 2456, 2569, 2579, 3456, 3468, 3478, 125a, 126a, 134b, 136b, 146b, 156a, 257a, 269a\\

$M_{68}$ &   1234, 1235, 1246, 1356, 2346, 2356, 2569, 2579, 3468, 3478, 125a, 126a, 134b, 136b, 146b, 156a, 257a, 269a\\
[1ex] % [1ex] adds vertical space
\hline %inserts single line
\end{tabular}
\label{table4} % is used to refer this table in the text
\end{table}

\section{Centrally symmetric quadrangulation of orientable surfaces}

A {\em cube} is a { 3}-dimensional solid object bounded by six square faces. It is one of the five {\em Platonic solids}. We consider boundary of the cube, that is, its boundary faces and denote it by {\em M}. Let $M :=\{$$[a_{000},$ $a_{001},$ $a_{101},$ $a_{100}],$ $[a_{000},$ $a_{100},$ $a_{110},$ $a_{010}],$ $[a_{000},$ $a_{001},$ $a_{011},$ $ a_{010}],$ $[a_{111},$ $a_{101},$ $a_{100},$ $a_{110}],$ $[a_{111},$ $a_{011},$ $a_{010},$ $a_{110}],$ $[a_{111},$ $a_{101},$ $a_{001},$ $a_{011}]\}$. The map {\em M} is a quadrangulation of sphere. Also, it is CS under the involution $I_{M} := a_{xyz} \mapsto a_{(1-x)(1-y(1-z)}$ where $x,y,z \in \{0,1\}$. Let $[a_{x_{1}x_{2}x_{3}},$ $a_{y_{1}y_{2}y_{3}},$ $a_{z_{1}z_{2}z_{3}},$ $a_{w_{1}w_{2}w_{3}}]$ be a face of $M_{0}$, see Figure 2.

%\bigskip

%\vspace{60 mm}
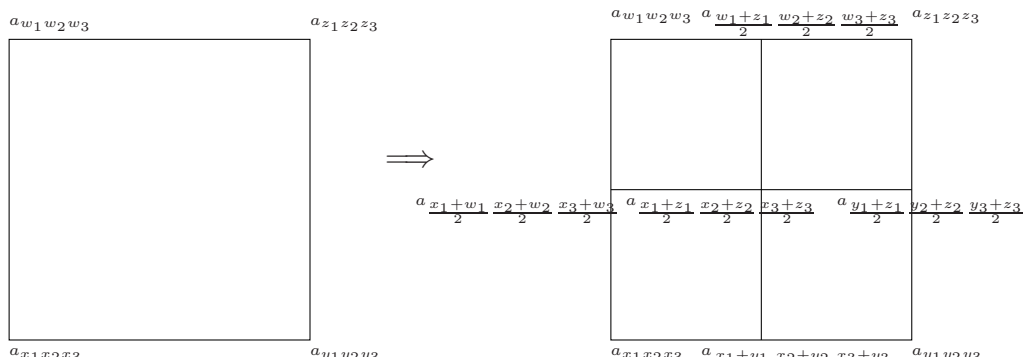
\begin{center}
\begin{picture}(50,0)(5,60)
\setlength{\unitlength}{10mm}
\drawpolygon (-4,2)(0,2)(0,6)(-4,6)
\put(-4,1.8) {\tiny $a_{x_{1}x_{2}x_{3}}$}
\put(0,1.8) {\tiny $a_{y_{1}y_{2}y_{3}}$}
\put(0,6.2) {\tiny $a_{z_{1}z_{2}z_{3}}$}
\put(-4,6.2) {\tiny $a_{w_{1}w_{2}w_{3}}$}

\put(-4,1){\tiny $[a_{x_{1}x_{2}x_{3}},$ $a_{y_{1}y_{2}y_{3}},$ $a_{z_{1}z_{2}z_{3}},$ $a_{w_{1}w_{2}w_{3}}]$}

%\drawpolygon (.5,0)(4.5,0)(2.75,2)
\put(1,4.3){$\Longrightarrow$}

\drawpolygon (4,2)(8,2)(8,6)(4,6)
\put(4,1.8) {\tiny $a_{x_{1}x_{2}x_{3}}$}
\put(8,1.8) {\tiny $a_{y_{1}y_{2}y_{3}}$}
\put(8,6.3) {\tiny $a_{z_{1}z_{2}z_{3}}$}
\put(4,6.3) {\tiny $a_{w_{1}w_{2}w_{3}}$}

\drawline[AHnb=0](4,4)(8,4)
\drawline[AHnb=0](6,2)(6,6)

\put(1.4,3.8) {\tiny $a_{\frac{x_{1}+w_{1}}{2}\frac{x_{2}+w_{2}}{2}\frac{x_{3}+w_{3}}{2}}$}
\put(4.2,3.8) {\tiny $a_{\frac{x_{1}+z_{1}}{2}\frac{x_{2}+z_{2}}{2}\frac{x_{3}+z_{3}}{2}}$}
\put(7,3.8) {\tiny $a_{\frac{y_{1}+z_{1}}{2}\frac{y_{2}+z_{2}}{2}\frac{y_{3}+z_{3}}{2}}$}
\put(5.2,6.3) {\tiny $a_{\frac{w_{1}+z_{1}}{2}\frac{w_{2}+z_{2}}{2}\frac{w_{3}+z_{3}}{2}}$}
\put(5.2,1.8) {\tiny $a_{\frac{x_{1}+y_{1}}{2}\frac{x_{2}+y_{2}}{2}\frac{x_{3}+y_{3}}{2}}$}

\put(3, 0){\tiny { Figure 2}}

\end{picture}
\end{center}
\vspace{-1in}
%\hrule
\vspace{4in}
%\newpage
We introduce a vertex $a_{\frac{x_{1}+y_{1}}{2}\frac{x_{2}+y_{2}}{2}\frac{x_{3}+y_{3}}{2}}$ which is middle vertex of the edge $a_{x_{1}x_{2}x_{3}}$ $a_{y_{1}y_{2}y_{3}}$. We introduce five vertices in each face and divide the face into four  4-gonal faces, see  Figure 2. So, we get a quadrangulated map, say $M_{0}$, from {\em M} where $M_{0} : = \{ [a_{000},$ $a_{00\frac{1}{2}},$ $ a_{\frac{1}{2}0\frac{1}{2}},$ $ a_{\frac{1}{2}00}],$ $[a_{00\frac{1}{2}},$ $a_{\frac{1}{2}0\frac{1}{2}},$ $ a_{\frac{1}{2}01},$ $a_{001}],$ $[a_{\frac{1}{2}01},$ $a_{\frac{1}{2}0\frac{1}{2}},$ $a_{10\frac{1}{2}},$ $a_{101}],$ $[a_{10\frac{1}{2}},$ $a_{\frac{1}{2}0\frac{1}{2}},$ $a_{\frac{1}{2}00}, a_{100}],$ $[a_{000},$ $a_{\frac{1}{2}00},$ $a_{\frac{1}{2}\frac{1}{2}0},$ $a_{0\frac{1}{2}0}],$ $[a_{100}, a_{\frac{1}{2}00},$ $a_{\frac{1}{2}\frac{1}{2}0},$ $ a_{1\frac{1}{2}0}],$ $[a_{110},$ $a_{1\frac{1}{2}0},$ $a_{\frac{1}{2}\frac{1}{2}0},$ $a_{\frac{1}{2}10}],$ $[a_{010},$ $ a_{0\frac{1}{2}0},$ $a_{\frac{1}{2}\frac{1}{2}0},$ $a_{\frac{1}{2}10}],$ $[a_{000},$ $a_{00\frac{1}{2}},$ $ a_{0\frac{1}{2}\frac{1}{2}},$ $a_{0\frac{1}{2}0}],$ $[a_{0\frac{1}{2}0},$ $a_{0\frac{1}{2}\frac{1}{2}},$ $a_{01\frac{1}{2}}, a_{010}]$ $[a_{011},$ $a_{0\frac{1}{2}1},$ $a_{0\frac{1}{2}\frac{1}{2}},$ $a_{01\frac{1}{2}}],$ $[a_{0\frac{1}{2}1},$ $ a_{0\frac{1}{2}\frac{1}{2}},$ $a_{00\frac{1}{2}},$ $a_{001}],$ $[a_{111},$ $a_{11\frac{1}{2}},$ $a_{\frac{1}{2}1\frac{1}{2}},$ $a_{\frac{1}{2}11}],$ $[a_{11\frac{1}{2}},$ $a_{\frac{1}{2}1\frac{1}{2}},$ $a_{\frac{1}{2}10},$ $a_{110}],$ $[a_{\frac{1}{2}10},$ $ a_{\frac{1}{2}1\frac{1}{2}},$ $a_{01\frac{1}{2}},$ $a_{010}],$ $[a_{01\frac{1}{2}},$ $a_{\frac{1}{2}1\frac{1}{2}},$ $a_{\frac{1}{2}11},$ $a_{011}],$ $[a_{111},$ $a_{\frac{1}{2}11},$ $a_{\frac{1}{2}\frac{1}{2}1},$ $a_{1\frac{1}{2}1}],$ $[a_{011},$ $ a_{\frac{1}{2}11},$ $a_{\frac{1}{2}\frac{1}{2}1},$ $a_{0\frac{1}{2}1}],$ $[a_{001},$ $a_{0\frac{1}{2}1},$ $a_{\frac{1}{2}\frac{1}{2}1},$ $a_{\frac{1}{2}01}],$ $[a_{101},$ $a_{1\frac{1}{2}1},$ $a_{\frac{1}{2}\frac{1}{2}1},$ $ a_{\frac{1}{2}01}],$ $[a_{111},$ $a_{11\frac{1}{2}},$ $a_{1\frac{1}{2}\frac{1}{2}},$ $a_{1\frac{1}{2}1}],$ $[a_{1\frac{1}{2}1},$ $ a_{1\frac{1}{2}\frac{1}{2}},$ $a_{10\frac{1}{2}}, a_{101}],$ $[a_{100},$ $a_{1\frac{1}{2}0},$ $a_{1\frac{1}{2}\frac{1}{2}},$ $ a_{10\frac{1}{2}}],$ $[a_{1\frac{1}{2}0},$ $a_{1\frac{1}{2}\frac{1}{2}},$ $a_{11\frac{1}{2}},$ $a_{110}]\}$. The map $M_{0}$ is centrally symmetric under the involution $I_{M_{0}} := a_{xyz}\mapsto a_{(1-x)(1-y(1-z)}$ where $x,y,z \in \{0,\frac{1}{2},1\}$ on the set of vertices $\{a_{000},$ $a_{001},$ $a_{101},$ $a_{100},$ $a_{011},$ $a_{010},$ $a_{110},$ $a_{111},$ $a_{00\frac{1}{2}},$ $a_{\frac{1}{2}0\frac{1}{2}},$ $a_{\frac{1}{2}00},$ $a_{\frac{1}{2}01},$ $a_{10\frac{1}{2}},$ $a_{0\frac{1}{2}0},$ $a_{\frac{1}{2}\frac{1}{2}0},$ $a_{1\frac{1}{2}0},$ $a_{\frac{1}{2}10},$ $a_{01\frac{1}{2}}$ $a_{0\frac{1}{2}\frac{1}{2}},$ $a_{11\frac{1}{2}},$ $a_{\frac{1}{2}1\frac{1}{2}},$ $a_{\frac{1}{2}11},$ $a_{1\frac{1}{2}1},$ $a_{\frac{1}{2}\frac{1}{2}1},$ $a_{1\frac{1}{2}1},$ $a_{1\frac{1}{2}\frac{1}{2}} \}$. We consider an isomorphic copy of $M_{0}$ and denote it by $N_{0}$. Let $N_{0} : = \{ [A_{000},$ $A_{00\frac{1}{2}},$ $A_{\frac{1}{2}0\frac{1}{2}},$ $A_{\frac{1}{2}00}],$ $[A_{00\frac{1}{2}},$ $ A_{\frac{1}{2}0\frac{1}{2}},$ $A_{\frac{1}{2}01},$ $A_{001}],$ $[A_{\frac{1}{2}01},$ $A_{\frac{1}{2}0\frac{1}{2}},$ $ A_{10\frac{1}{2}},$ $A_{101}],$ $[A_{10\frac{1}{2}},$ $A_{\frac{1}{2}0\frac{1}{2}},$ $A_{\frac{1}{2}00},$ $A_{100}],$ $[A_{000},$ $A_{\frac{1}{2}00},$ $A_{\frac{1}{2}\frac{1}{2}0},$ $ A_{0\frac{1}{2}0}],$ $[A_{100},$ $A_{\frac{1}{2}00},$ $A_{\frac{1}{2}\frac{1}{2}0},$ $A_{1\frac{1}{2}0}],$ $[A_{110},$ $A_{1\frac{1}{2}0},$ $A_{\frac{1}{2}\frac{1}{2}0},$ $ A_{\frac{1}{2}10}],$ $[A_{010},$ $A_{0\frac{1}{2}0},$ $A_{\frac{1}{2}\frac{1}{2}0},$ $A_{\frac{1}{2}10}],$ $[A_{000}, A_{00\frac{1}{2}},$ $A_{0\frac{1}{2}\frac{1}{2}},$ $A_{0\frac{1}{2}0}],$ $[A_{0\frac{1}{2}0},$ $A_{0\frac{1}{2}\frac{1}{2}},$ $A_{01\frac{1}{2}},$ $A_{010}],$ $[A_{011},$ $A_{0\frac{1}{2}1},$ $A_{0\frac{1}{2}\frac{1}{2}},$ $A_{01\frac{1}{2}}]$ $[A_{0\frac{1}{2}1},$ $ A_{0\frac{1}{2}\frac{1}{2}},$ $A_{00\frac{1}{2}},$ $A_{001}]$ $[A_{111},$ $A_{11\frac{1}{2}},$ $ A_{\frac{1}{2}1\frac{1}{2}},$ $A_{\frac{1}{2}11}],$ $[A_{11\frac{1}{2}},$ $A_{\frac{1}{2}1\frac{1}{2}},$ $A_{\frac{1}{2}10},$ $A_{110}],$ $[A_{\frac{1}{2}10},$ $A_{\frac{1}{2}1\frac{1}{2}},$ $A_{01\frac{1}{2}},$ $A_{010}],$ $[A_{01\frac{1}{2}},$ $ A_{\frac{1}{2}1\frac{1}{2}},$ $A_{\frac{1}{2}11},$ $A_{011}],$ $[A_{111},$ $A_{\frac{1}{2}11},$ $A_{\frac{1}{2}\frac{1}{2}1},$ $ A_{1\frac{1}{2}1}],$ $[A_{011},$ $A_{\frac{1}{2}11},$ $A_{\frac{1}{2}\frac{1}{2}1},$ $A_{0\frac{1}{2}1}],$ $[A_{001},$ $A_{0\frac{1}{2}1},$ $A_{\frac{1}{2}\frac{1}{2}1},$ $A_{\frac{1}{2}01}],$ $[A_{101},$ $A_{1\frac{1}{2}1},$ $ A_{\frac{1}{2}\frac{1}{2}1},$ $A_{\frac{1}{2}01}],$ $[A_{111},$ $A_{11\frac{1}{2}},$ $A_{1\frac{1}{2}\frac{1}{2}},$ $ A_{1\frac{1}{2}1}],$ $[A_{1\frac{1}{2}1},$ $ A_{1\frac{1}{2}\frac{1}{2}},$ $A_{10\frac{1}{2}},$ $A_{101}]$ $[A_{100}, A_{1\frac{1}{2}0},$ $A_{1\frac{1}{2}\frac{1}{2}},$ $A_{10\frac{1}{2}}],$ $[A_{1\frac{1}{2}0}, A_{1\frac{1}{2}\frac{1}{2}}, $ $A_{11\frac{1}{2}},$ $A_{110}]\}$. It is CS under the involution $I_{N_{0}} := A_{xyz}\mapsto A_{(1-x)(1-y(1-z)}$ where $x,y,z \in \{0,\frac{1}{2},1\}$ on the set of vertices $\{A_{000},$ $A_{001},$ $A_{101},$ $A_{100},$ $A_{011},$ $A_{010},$ $A_{110},$ $A_{111},$ $A_{00\frac{1}{2}},$ $A_{\frac{1}{2}0\frac{1}{2}},$ $A_{\frac{1}{2}00},$ $A_{\frac{1}{2}01},$ $A_{10\frac{1}{2}},$ $A_{0\frac{1}{2}0},$ $A_{\frac{1}{2}\frac{1}{2}0},$ $A_{1\frac{1}{2}0},$ $A_{\frac{1}{2}10},$ $A_{01\frac{1}{2}}$ $A_{0\frac{1}{2}\frac{1}{2}},$ $A_{11\frac{1}{2}},$ $A_{\frac{1}{2}1\frac{1}{2}},$ $A_{\frac{1}{2}11},$ $A_{1\frac{1}{2}1},$ $A_{\frac{1}{2}\frac{1}{2}1},$ $A_{1\frac{1}{2}1},$ $A_{1\frac{1}{2}\frac{1}{2}} \}$. We remove interior of two faces namely $[a_{000}, a_{00\frac{1}{2}}, a_{\frac{1}{2}0\frac{1}{2}}, a_{\frac{1}{2}00}],$ $[a_{111}, a_{11\frac{1}{2}}, a_{\frac{1}{2}1\frac{1}{2}}, a_{\frac{1}{2}11}]$ from $M_{0}$ and $[A_{000},$ $ A_{00\frac{1}{2}},$ $A_{\frac{1}{2}0\frac{1}{2}},$ $A_{\frac{1}{2}00}],$ $[A_{111},$ $A_{11\frac{1}{2}},$ $A_{\frac{1}{2}1\frac{1}{2}},$ $ A_{\frac{1}{2}11}]$ from $N_{0}$. We identify the cycles $C(a_{000},$ $a_{00\frac{1}{2}},$ $a_{\frac{1}{2}0\frac{1}{2}},$ $a_{\frac{1}{2}00})$ with $C(A_{000},$ $A_{00\frac{1}{2}},$ $A_{\frac{1}{2}0\frac{1}{2}},$ $A_{\frac{1}{2}00})$ and $C(a_{111},$ $a_{11\frac{1}{2}},$ $a_{\frac{1}{2}1\frac{1}{2}},$ $a_{\frac{1}{2}11})$ with $C(A_{111},$ $A_{11\frac{1}{2}},$ $A_{\frac{1}{2}1\frac{1}{2}},$ $A_{\frac{1}{2}11})$ by the map $a_{xyz}\mapsto A_{xyz}$, where $xyz \in \{000,$ $00\frac{1}{2},$ $\frac{1}{2}0\frac{1}{2},$ $\frac{1}{2}00,$ $111,$ $11\frac{1}{2},$ $\frac{1}{2}1\frac{1}{2},$ $\frac{1}{2}11\}$. Here we identify eight vertices $a_{000},$ $a_{00\frac{1}{2}},$ $ a_{\frac{1}{2}0\frac{1}{2}},$ $a_{\frac{1}{2}00}$ $a_{111},$ $a_{11\frac{1}{2}},$ $a_{\frac{1}{2}1\frac{1}{2}},$ $a_{\frac{1}{2}11}$ with the eight vertices $A_{000},$ $A_{00\frac{1}{2}},$ $A_{\frac{1}{2}0\frac{1}{2}},$ $A_{\frac{1}{2}00}$ $A_{111},$ $A_{11\frac{1}{2}},$ $A_{\frac{1}{2}1\frac{1}{2}},$ $A_{\frac{1}{2}11}$ respectively. Hence we get $\overline{M}_{1} := \{ [A_{00\frac{1}{2}},$ $A_{\frac{1}{2}0\frac{1}{2}},$ $ a_{\frac{1}{2}01},$ $a_{001}],$ $[a_{\frac{1}{2}01},$ $A_{\frac{1}{2}0\frac{1}{2}},$ $a_{10\frac{1}{2}}, a_{101}],$ $[a_{10\frac{1}{2}},$ $ A_{\frac{1}{2}0\frac{1}{2}},$ $A_{\frac{1}{2}00},$ $a_{100}],$ $[A_{000}, A_{\frac{1}{2}00}, a_{\frac{1}{2}\frac{1}{2}0},$ $a_{0\frac{1}{2}0}],$ $[a_{100},$ $A_{\frac{1}{2}00},$ $a_{\frac{1}{2}\frac{1}{2}0},$ $ a_{1\frac{1}{2}0}],$ $[a_{110},$ $a_{1\frac{1}{2}0},$ $a_{\frac{1}{2}\frac{1}{2}0},$ $a_{\frac{1}{2}10}],$ $[a_{010},$ $a_{0\frac{1}{2}0},$ $a_{\frac{1}{2}\frac{1}{2}0},$ $a_{\frac{1}{2}10}],$ $[A_{000},$ $A_{00\frac{1}{2}},$ $a_{0\frac{1}{2}\frac{1}{2}},$ $a_{0\frac{1}{2}0}],$ $[a_{0\frac{1}{2}0},$ $a_{0\frac{1}{2}\frac{1}{2}},$ $a_{01\frac{1}{2}}, a_{010}],$ $[a_{011},$ $a_{0\frac{1}{2}1},$ $a_{0\frac{1}{2}\frac{1}{2}},$ $a_{01\frac{1}{2}}]$ $[a_{0\frac{1}{2}1}, a_{0\frac{1}{2}\frac{1}{2}},$ $A_{00\frac{1}{2}},$ $a_{001}],$ $[A_{11\frac{1}{2}},$ $A_{\frac{1}{2}1\frac{1}{2}},$ $a_{\frac{1}{2}10},$ $a_{110}],$ $[a_{\frac{1}{2}10},$ $A_{\frac{1}{2}1\frac{1}{2}},$ $a_{01\frac{1}{2}},$ $a_{010}],$ $[a_{01\frac{1}{2}},$ $ A_{\frac{1}{2}1\frac{1}{2}},$ $A_{\frac{1}{2}11},$ $a_{011}],$ $[A_{111}, A_{\frac{1}{2}11}, a_{\frac{1}{2}\frac{1}{2}1},$ $a_{1\frac{1}{2}1}],$ $[a_{011},$ $A_{\frac{1}{2}11},$ $a_{\frac{1}{2}\frac{1}{2}1},$ $ a_{0\frac{1}{2}1}],$ $[a_{001},$ $a_{0\frac{1}{2}1},$ $a_{\frac{1}{2}\frac{1}{2}1},$ $a_{\frac{1}{2}01}],$ $[a_{101},$ $a_{1\frac{1}{2}1},$ $a_{\frac{1}{2}\frac{1}{2}1},$ $a_{\frac{1}{2}01}],$ $[A_{111},$ $A_{11\frac{1}{2}},$ $ a_{1\frac{1}{2}\frac{1}{2}},$ $a_{1\frac{1}{2}1}],$ $[a_{1\frac{1}{2}1},$ $a_{1\frac{1}{2}\frac{1}{2}},$ $a_{10\frac{1}{2}},$ $ a_{101}],$ $[a_{100},$ $a_{1\frac{1}{2}0},$ $a_{1\frac{1}{2}\frac{1}{2}},$ $a_{10\frac{1}{2}}],$ $[a_{1\frac{1}{2}0},$ $ a_{1\frac{1}{2}\frac{1}{2}},$ $A_{11\frac{1}{2}},$ $a_{110}],$ $[A_{00\frac{1}{2}},$ $A_{\frac{1}{2}0\frac{1}{2}},$ $A_{\frac{1}{2}01},$ $A_{001}],$ $[A_{\frac{1}{2}01},$ $A_{\frac{1}{2}0\frac{1}{2}},$ $A_{10\frac{1}{2}},$ $A_{101}],$ $[A_{10\frac{1}{2}},$ $A_{\frac{1}{2}0\frac{1}{2}},$ $A_{\frac{1}{2}00},$ $A_{100}],$ $[A_{000},$ $A_{\frac{1}{2}00},$ $ A_{\frac{1}{2}\frac{1}{2}0},$ $A_{0\frac{1}{2}0}],$ $[A_{100}, A_{\frac{1}{2}00},$ $A_{\frac{1}{2}\frac{1}{2}0},$ $A_{1\frac{1}{2}0}],$ $[A_{110}, A_{1\frac{1}{2}0},$ $A_{\frac{1}{2}\frac{1}{2}0},$ $A_{\frac{1}{2}10}],$ $[A_{010},$ $A_{0\frac{1}{2}0},$ $A_{\frac{1}{2}\frac{1}{2}0},$ $A_{\frac{1}{2}10}],$ $[A_{000},$ $A_{00\frac{1}{2}},$ $ A_{0\frac{1}{2}\frac{1}{2}},$ $A_{0\frac{1}{2}0}],$ $[A_{0\frac{1}{2}0},$ $A_{0\frac{1}{2}\frac{1}{2}},$ $A_{01\frac{1}{2}},$ $ A_{010}],$ $[A_{011},$ $A_{0\frac{1}{2}1},$ $A_{0\frac{1}{2}\frac{1}{2}},$ $A_{01\frac{1}{2}}],$ $[A_{0\frac{1}{2}1}, A_{0\frac{1}{2}\frac{1}{2}},$ $A_{00\frac{1}{2}},$ $A_{001}],$ $[A_{11\frac{1}{2}},$ $A_{\frac{1}{2}1\frac{1}{2}},$ $A_{\frac{1}{2}10},$ $A_{110}],$ $[A_{\frac{1}{2}10}, A_{\frac{1}{2}1\frac{1}{2}},$ $A_{01\frac{1}{2}},$ $A_{010}],$ $[A_{01\frac{1}{2}},$ $A_{\frac{1}{2}1\frac{1}{2}},$ $A_{\frac{1}{2}11},$ $A_{011}],$ $[A_{111},$ $A_{\frac{1}{2}11},$ $ A_{\frac{1}{2}\frac{1}{2}1},$ $A_{1\frac{1}{2}1}],$ $[A_{011},$ $A_{\frac{1}{2}11},$ $ A_{\frac{1}{2}\frac{1}{2}1},$ $A_{0\frac{1}{2}1}],$ $[A_{001},$ $A_{0\frac{1}{2}1},$ $A_{\frac{1}{2}\frac{1}{2}1},$ $A_{\frac{1}{2}01}],$ $[A_{101}, A_{1\frac{1}{2}1},$ $A_{\frac{1}{2}\frac{1}{2}1},$ $A_{\frac{1}{2}01}],$ $[A_{111},$ $A_{11\frac{1}{2}},$ $A_{1\frac{1}{2}\frac{1}{2}},$ $A_{1\frac{1}{2}1}],$ $[A_{1\frac{1}{2}1},$ $A_{1\frac{1}{2}\frac{1}{2}},$ $A_{10\frac{1}{2}},$ $A_{101}]$ $[A_{100},$ $A_{1\frac{1}{2}0},$ $A_{1\frac{1}{2}\frac{1}{2}},$ $A_{10\frac{1}{2}}],$ $[A_{1\frac{1}{2}0},$ $A_{1\frac{1}{2}\frac{1}{2}},$ $A_{11\frac{1}{2}},$ $A_{110}]\}$. We define $I_{\overline{M}_{1}} :=I_{M_{0}}I_{N_{0}}-\{(a_{000},$ $a_{111}),$ $(a_{00\frac{1}{2}},$ $a_{11\frac{1}{2}}),$ $(a_{\frac{1}{2}0\frac{1}{2}},$  $a_{\frac{1}{2}1\frac{1}{2}}),$ $(a_{\frac{1}{2}00},$ $a_{\frac{1}{2}11})\}$.
The map $\overline{M}_{1}$ is CS under the involution $I_{\overline{M}_{1}}$. The Euler characteristic of $\overline{M}_{1}$ is $\chi(M_{0})+\chi(N_{0})-4 = 2+2-4 = 0$. Therefore the map $\overline{M}_{1} := M_{0} \# N_{0}$ ( this is usual notation of connected sum ) is centrally symmetric quadrangulation of torus on 1.(26 - 8)+ 26 = 1.18+26 = 44 vertices.

 Next, we denote $M_{1} := \overline{M}_{1}$ and $I_{M_{1}}:= I_{\overline{M}_{1}}$. We consider two faces $F_{1}$, $F_{2}$ of $M_{1}$ where $F_{2}^{I_{M_{1}}} = F_{1}$. Also, we consider an isomorphic copy of $M_{0}$ and denote it by $N_{1}$. The map $N_{1}$ is CS under $I_{N_{1}}:= I_{M_{0}}$. We choose two faces $F_{3}, F_{4}$ of $N_{1}$ where $F_{3}^{I_{N_{1}}} = F_{4}$. We remove interior of $F_{1}$ and $F_{2}$ from $M_{1}$ and $F_{3}$ and $F_{4}$ from $N_{1}$. Similarly, we define $I_{\overline{M}_{1}}$ using $I_{M_{1}}$ and $I_{N_{1}}$. Also we identify $\partial F_{1}$ with $\partial F_{3}$ and $\partial F_{2}$ with $\partial F_{4}$. Hence we get an object $\overline{M}_{2}$ which is CS under $I_{\overline{M}_{2}}$. Where the map $\overline{M}_{2} := M_{1}\#N_{1}$ of orientable genus $g = 2$ on $2.(26-8)+26 = 2.18+26 = 62$ vertices.

 Similarly, at $g^{th}$ step, we consider the object which is obtained at $(g-1)^{th}$ step and denote it by $M_{g-1}$. We also denote the involution of the object under which the object is CS and denote it by $I_{M_{g-1}}$. We consider an isomorphic copy of $M_{0}$ and denote it by $N_{g-1}$. It is CS under $I_{N_{g-1}}:= I_{1}$. By the similar process, we get a map $\overline{M}_{g} := M_{g-1}\#N_{g-1}$ on $g.(26-8)+26 = g.18+26$ vertices. Similarly, we define $I_{\overline{M}_{g}}$ using $I_{M_{g-1}}$ and $I_{N_{g-1}}$. The map $\overline{M}_{g}$ is CS under the involution $I_{\overline{M}_{g}}$. Therefore, we can construct centrally symmetric orientable quadrangulated surfaces for any positive orientable genus.

\bigskip
\noindent{\sc Proof of Theorem}\ref{thm2} The proof of theorem \ref{thm2} now follows from arguments in this section.
\hfill$\Box$

\section{Centrally symmetric orientable surfaces all of whose faces are pentagons}

A {\em dodecahedron} is a {\em polyhedron} with twelve faces. It is one of the five {\em Platonic solids}. Its boundary composite of { 12} regular pentagonal faces with three meeting at each vertex. We consider the boundary of the dodecahedron and denote it by $M_{0}$. Let $M_{0} :=$ $\{[1,$ $2,$ $17,$ $16,$ $10],$ $[2,$ $3,$ $4,$ $18,$ $ 17],$ $[4,$ $5,$ $6,$ $19,$ $18],$ $[6,$ $19,$ $20,$ $8,$ $7],$ $[8,$ $20,$ $16,$ $10,$ $9],$ $[16,$ $17,$ $18,$ $19,$ $20],$ $[1,$ $2, 3,$ $12,$ $11],$ $[3,$ $4,$ $5,$ $13,$ $12],$ $[5,$ $6,$ $7,$ $14,$ $13],$ $[7, 8, 9, 15, 14],$ $[9,$ $10,$ $1,$ 11,$ 15],$ $[11, 12, 13, 14, 15]\}$. The map $M_{0}$ is CS under the involution $I_{M_{0}} =$ $(5,$ $10)$ $(2,$ $7)$ $(3,$ $8)$ $(1,$ $6)$ $(4,$ $9)$ $(14,$ $17)$ $(13,$ $16)$ $(15,$ $18)$ $(11,$ $19)$ $(12,$ $20)$. We consider an isomorphic copy of $M_{0}$ and denote it by $N_{0}$. Let $N_{0} :=$ $\{[21,$ $22,$ $37,$ $36,$ $30],$ $[22,$ $23,$ $24,$ $38,$ $37],$ $[24,$ $25,$ $26,$ $39,$ $38],$ $[26,$ $39,$ $40,$ $ 28,$ $27],$ $[28,$ $40,$ $36,$ $30,$ $29],$ $[36,$ $37,$ $38,$ $39,$ $40],$ $[21,$ $22,$ $23,$ $32,$ $31],$ $[23,$ $24,$ $25,$ $33,$ $32],$ $[25,$ $26,$ $27,$ $34,$ $33],$ $[27,$ $28,$ $29,$ $35,$ $34],$ $[29,$ $30,$ $21,$ $31,$ $35],$ $[31,$ $32,$ $33,$ 34,$ $35]$ \}$. It is CS under the involution $I_{N_{0}} =$ $(25,$ $30)$ $(22,$ $27)$ $(23,$ $28)$ $(21,$ $26)$ $(24,$ $29)$ $(34,$ $37)$ $(33,$ $36)$ $(35,$ $38)$ $(31,$ $39)$ $(32,$ $40)$. We remove interior of two faces namely $[1,$ $2,$ $17,$ $16,$ $10],$ $[6,$ $7,$ $14,$ $13,$ $5]$ from $M_{0}$ and $[21,$ $22,$ $37,$ $36,$ $30],$ $[26,$ $27,$ $34,$ $33,$ $25]$ from $N_{0}$. We identify the cycles $C(1,$ $2,$ $17,$ $16,$ $10)$ with $C(21,$ $22,$ $37,$ $36,$ $30)$ by the map $1 \mapsto 21,$ $2 \mapsto 22,$ $17 \mapsto 37,$ $16 \mapsto 36,$ $10 \mapsto 30$ and $C(6, 7, 14, 13, 5)$ with $C(26, 27, 34, 33, 25)$ by the map $6 \mapsto 26,$ $7 \mapsto 27,$ $14 \mapsto 34,$ $13 \mapsto 33,$ $5 \mapsto 25$. In this process we identify ten vertices $1,$ $2,$ $17,$ $16,$ $10,$ $6,$ $7,$ $14,$ $13,$ $5$ with ten vertices $21,$ $22,$ $37,$ $36,$ $30,$ $26,$ $27,$ $34,$ $33,$ $25$. Hence we get an object and denote it by $\overline{M}_{1}$ where $\overline{M}_{1} :=$ $\{ [22,$ $3,$ $4,$ $18,$ $37],$ $[4,$ $25,$ $26,$ $19,$ $18],$ $[26,$ $19,$ $20,$ $8,$ $27],$ $[8,$ $20,$ $36,$ $30,$ $9],$ $[36,$ $37,$ $18,$ $19,$ $20],$ $[21,$ $22,$ $3,$ $12,$ $11],$ $[3,$ $4,$ $25,$ $33,$ $12],$ $[27,$ $8,$ $9,$ $15,$ $34],$ $[9,$ $30,$ $21,$ $11,$ $15],$ $[11,$ $12,$ $33,$ $34,$ $15],$ $[22,$ $23,$ $24,$ $38,$ $37],$ $[24,$ $25,$ $26,$ $39,$ $38],$ $[26,$ $39,$ $40,$ $28,$ $27],$ $[28,$ $40,$ $36,$ $30,$ $29],$ $[36,$ $37,$ $38,$ $39,$ $40],$ $[21,$ $22,$ 23,$ $32,$ 31],$ $[23,$ $24,$ $25,$ $33,$ $32],$ $[27,$ $28,$ $29,$ $35,$ $34],$ $[29,$ $30,$ $21,$ $31,$ $35],$ $[31,$ $32,$ $33,$ $ 34,$ $35] \}$. We define $I_{\overline{M}_{1}} :=$ $(3,$ $8)$ $(4,$ $9)$ $(15,$ $18)$ $(11,$ $19)$ $(12,$ $20)$ $(25,$ $30)$ $(22,$ $27)$ $(23,$ $28)$ $(21,$ $26)$ $(24,$ $29)$ $(34,$ $37)$ $(33,$ $36)$ $(35,$ $38)$ $(31,$ $39)$ $(32,$ $40)$. Therefore the map $\overline{M}_{1} := M_{0}\#N_{0}$ is CS under the involution $I_{\overline{M}_{1}}$ on $1.(20-10)+20 = 1.10+20 = 30$ vertices. The Euler characteristic of $\overline{M}_{1}$ is $\chi(M_{0})+\chi(N_{0})-4 = 2+2-4 = 0$. Next, we follow the { Section 3}. Hence we get a CS map $\overline{M}_{2}$ where $\overline{M}_{2} := M_{1}\#N_{1}$ on $2.(20-10)+20 = 2.10+20 = 40$ vertices. Similarly, at $g^{th}$ step, we get CS map $\overline{M}_{g} := M_{g-1}\#N_{g-1}$ of genus $g$ on $g.(20-10)+20 = g.10+20$ vertices. Therefore, we can construct CS orientable surfaces all of whose faces are pentagons for any positive genus.

\bigskip

\noindent{\sc Proof of Theorem}\ref{thm3} The proof of theorem \ref{thm3} now follows from arguments in this section.
\hfill$\Box$

\section{Centrally symmetric orientable surfaces all of whose faces are hexagons}

Let $M_{1}$ be a CS torus all of whose faces are hexagons on { 24} vertices. Let $M_{1} :=$ $\{ [1,$ $2,$ $3,$ $8,$ $7,$ $6],$ $[3,$ $4,$ $5,$ $10,$ $9,$ $8],$ $[5,$ $6,$ $7,$ $12,$ $11,$ $10],$ $[7,$ $8,$ $9,$ $14,$ $13,$ $12],$ $[9,$ $10,$ $11,$ $16,$ $15,$ $14],$ $[11,$ $12,$ $13,$ $18,$ $17,$ $16],$ $[13,$ $14,$ $15,$ $20,$ $19,$ $18],$ $[15,$ $16,$ $17,$ $22,$ $21,$ $20],$ $[17,$ $18,$ $19,$ $24,$ $23,$ $22],$ $[19,$ $20,$ $21,$ $2,$ $1,$ $24],$ $[21,$ $22,$ $23,$ $4,$ $3,$ $2],$ $[23,$ $24,$ $1,$ $6,$ $5,$ $4] \}$. The map $M_{1}$ is CS under the involution $I_{M_{1}} =$ $(1,$ $13)$ $(2,$ $14)$ $(3,$ $15)$ $(4,$ $16)$ $(5,$ $17)$ $(6,$ $18)$ $(7,$ $19)$ $(8,$ $20)$ $(9,$ $21)$ $(10,$ $22)$ $(11,$ $23)$ $(12,$ $24)$. We consider an isomorphic copy of $M_{1}$ and denote it by $N_{1}$. Let $N_{1} :=$ $\{ [a_{1},$ $a_{2},$ $a_{3},$ $a_{8},$ $a_{7},$ $a_{6}],$ $[a_{3},$ $a_{4},$ $a_{5},$ $a_{10},$ $a_{9},$ $a_{8}],$ $[a_{5},$ $a_{6},$ $a_{7},$ $a_{12},$ $a_{11},$ $a_{10}],$ $[a_{7},$ $a_{8},$ $a_{9},$ $a_{14},$ $a_{13},$ $a_{12}],$ $[a_{9},$ $ a_{10},$ $a_{11},$ $a_{16},$ $a_{15},$ $a_{14}],$ $[a_{11},$ $a_{12},$ $a_{13},$ $a_{18},$ $a_{17},$ $a_{16}],$ $[a_{13},$ $a_{14},$ $a_{15},$ $a_{20},$ $a_{19},$ $a_{18}],$ $[a_{15},$ $a_{16},$ $a_{17},$ $a_{22},$ $a_{21},$ $a_{20}],$ $[a_{17},$ $a_{18},$ $a_{19},$ $a_{24},$ $a_{23},$ $a_{22}],$ $[a_{19},$ $a_{20},$ $a_{21},$ $a_{2},$ $a_{1},$ $a_{24}],$ $[a_{21},$ $a_{22},$ $a_{23},$ $a_{4},$ $a_{3},$ $a_{2}],$ $[a_{23},$ $a_{24},$ $a_{1},$ $a_{6},$ $a_{5},$ $a_{4}] \}$. It is CS under the involution $I_{N_{1}} =$ $(a_{1},$ $ a_{13})$ $(a_{2},$ $a_{14})$ $(a_{3},$ $a_{15})$ $(a_{4},$ $a_{16})$ $(a_{5},$ $a_{17})$ $(a_{6},$ $a_{18})$ $(a_{7},$ $a_{19})$ $(a_{8},$ $a_{20})$ $(a_{9},$ $a_{21})$ $(a_{10},$ $a_{22})$ $(a_{11},$ $a_{23})$ $(a_{12},$ $a_{24})$. We remove interior of two faces namely $[1,$ $2,$ $3,$ $8,$ $7,$ $6],$ $[13,$ $14,$ $15,$ $20,$ $19,$ $18]$ from $M_{1}$ and $[a_{1},$ $a_{2},$ $a_{3},$ $a_{8},$ $a_{7},$ $a_{6}],$ $[a_{13},$ $a_{14},$ $a_{15},$ $a_{20},$ $a_{19},$ $a_{18}]$ from $N_{1}$. We identify the cycles $C(1, 2, 3, 8, 7, 6)$ with $C(a_{1}, a_{2}, a_{3}, a_{8}, a_{7}, a_{6})$ by the map $a_{1}\mapsto 1,$ $a_{2}\mapsto 2,$ $a_{3}\mapsto 3,$ $a_{8}\mapsto 8,$ $a_{7}\mapsto 7,$ $a_{6}\mapsto 6$ and $C(13, 14, 15, 20, 19, 18)$ with $C(a_{13}, a_{14}, a_{15}, a_{20}, a_{19}, a_{18})$ by the map $ a_{13}\mapsto 13,$ $a_{14}\mapsto 14,$ $a_{15}\mapsto 15,$ $a_{20}\mapsto 20,$ $a_{19}\mapsto 19,$ $a_{18}\mapsto 18$. Here we identify twelve vertices $1,$ $2,$ $3,$ $8,$ $7,$ $6,$ $13,$ $14,$ $15,$ $20,$ $19,$ $18$ with twelve vertices $a_{1},$ $a_{2},$ $a_{3},$ $a_{8},$ $a_{7},$ $a_{6},$ $a_{13},$ $a_{14},$ $a_{15},$ $a_{20},$ $a_{19},$ $a_{18}$. Therefore, we get the map $\overline{M}_{3} :=$ $\{[1,$ $2,$ $3,$ $8,$ $7,$ $6],$ $[3,$ $a_{4},$ $a_{5},$ $a_{10},$ $a_{9},$ $a_{8}],$ $[a_{5},$ $6,$ $7,$ $a_{12},$ $a_{11},$ $a_{10}],$ $[7,$ $8,$ $a_{9},$ $14,$ $13,$ $a_{12}],$ $[a_{9},$ $a_{10},$ $a_{11},$ $a_{16},$ $15,$ $14],$ $[a_{11},$ $a_{12},$ $13,$ $18,$ $a_{17},$ $a_{16}],$ $[13,$ $14,$ $15,$ $20,$ $19,$ $18],$ $[15,$ $a_{16},$ $a_{17},$ $a_{22},$ $a_{21},$ $20],$ $[a_{17},$ $18,$ $19,$ $a_{24},$ $a_{23},$ $a_{22}],$ $[19,$ $20,$ $a_{21},$ $2,$ $1,$ $a_{24}],$ $[a_{21},$ $a_{22},$ $a_{23},$ $a_{4},$ $3,$ $2],$ $[a_{23},$ $a_{24},$ $1,$ $6,$ $a_{5},$ $a_{4}],$ $[1,$ $2,$ $ 3,$ $48,$ $7,$ $6],$ $[3,$ $4,$ $5,$ $10,$ $9,$ $8],$ $[5,$ $6,$ $7,$ $12,$ $11,$ $10],$ $[7,$ $8,$ $9,$ $14,$ $13,$ $12],$ $[9,$ $10,$ $11,$ $16,$ $15,$ $14],$ $[11,$ $12,$ $13,$ $18,$ $17,$ $16],$ $[13,$ $14,$ $15,$ $20,$ $19,$ $18],$ $[15,$ $16,$ $17,$ $22,$ $21,$ $ 20],$ $[17,$ $18,$ $19,$ $24,$ $23,$ $22],$ $[19,$ $20,$ $21,$ $2,$ $1,$ $24],$ $[21,$ $22,$ $23,$ $4,$ $3,$ $2],$ $[23,$ $24,$ $1,$ $ 6,$ $5,$ $4] \}$. We define $I_{\overline{M}_{3}} :=$ $(1,$ $13)$ $(2,$ $14)$ $(3,$ $15)$ $(4,$ $16)$ $(5,$ $17)$ $(6,$ $18)$ $(7,$ $19)$ $(8,$ $20)$ $(9,$ $21)$ $(10,$ $22)$ $(11,$ $23)$ $(12,$ $24)$ $(a_{4},$ $a_{16})$ $(a_{5},$ $a_{17})$ $(a_{9},$ $a_{21})$ $(a_{10},$ $a_{22})$ $(a_{11},$ $a_{23})$ $(a_{12},$ $a_{24})$. Therefore the map $\overline{M}_{3} : = M_{1}\# N_{1}$ is CS under the involution $I_{\overline{M}_{3}}$ on $24 + 1.(24-12) = 24 + 1.12 =36$ vertices. The Euler characteristic  $\chi(\overline{M}_{3}) = -4$, that is, orientable genus of the map $\overline{M}_{3}$ is three. We use $g(\overline{M}_{3})$ to denote the number of orientable genus of the map $\overline{M}_{3}$. We follow the idea of { Section 3}. Hence we get CS map $\overline{M}_{5}: = M_{3}\#N_{3}$ of genus $g(\overline{M}_{5}) = g(M_{3}) + g(N_{3}) + 1 = 3 + 1 + 1 = 5$ on $24 + 2.(24-12) = 24+ 2.12 = 48$ vertices. Similarly, at $k^{th}$ step, we get a CS map all of whose faces are hexagons of genus $g = 2k-1$ on $24 + [\frac{2k-1}{2}]12$ vertices.  Therefore we can construct centrally symmetric orientable surfaces all of whose faces are hexagons for any odd number of genus.

\bigskip

\noindent{\sc Proof of Theorem}\ref{thm4} The proof of theorem \ref{thm4} now follows from arguments in this section.
\hfill$\Box$

\section{Construction of CS surfaces of type $\{q, p\}$ from known CS surfaces of type $\{p, q\}$}

The dual map $\overline{M}$ of {\em M} is a map on same surface as {\em M}. The map $\overline{M}$ has for its vertices the set of faces of {\em M} and two vertices of {\em M} are ends of an edge of {\em M} if the corresponding faces in {\em M} have an edge in common.

Let {\em M} be a centrally symmetric map on surfaces of type $\{p, q\}$ with respect to an involution $ I = (1, 2m)(2, 2m-1) \dots (m, m+1)$. Let {\em F} be a face of {\em M}. The $I$ orbit of {\em F} contains exactly two disjoint faces. Let $k:=|O|$ denote the number of orbits of {\em M}. Then, {\em M} contains { 2k} number of faces which is even. We take dual of {\em M} and denote it by $\overline{M}$ of type $\{q, p\}$. Let $F_{i}$ be a face in {\em M}. We denote $u_{i}$ be the dual vertex of $F_{i}$ in $\overline{M}$. Let $F_{2k-(i-1)}$ be the face which belongs to the orbit of $F_{i}$. We define an involution $\overline{I} : =$ $(u_{1}, u_{2k})$ $(u_{2}, u_{2k-1})$ \dots $(u_{i}, u_{2k-(i-1)})$ \dots $(u_{k}, u_{k+1})$. We claim that $\overline{M}$ is centrally symmetric under the involution $\overline{I}$. Suppose there is a face {\em F} which is fix under the involution $\overline{I}$, that is, $F^{\overline{I}} = F$. Let the dual faces of $F^{\overline{I}}$ and {\em F} be $X^{I}$ and {\em X} in {\em M} respectively. By the definition of duality $X^{I} = X$. Hence {\em X} is a fixed face in {\em M} under {\em I}. This shows {\em M} is not centrally symmetric, which is a contradiction to our assumption. Therefore the map $\overline{M}$ is centrally symmetric under $\overline{I}$.

\begin{example}\label{example2}
{\small {\bf Construction of CS sphere of type $\{3, 4\}$ from known CS sphere of type $\{4, 3\}$ }}
\end{example}

We consider $M$ the boundary of the cube which is a map of type $\{4, 3\}$ on sphere. Let $M :=\{$ $[a_{000},$ $a_{001},$ $a_{101},$ $a_{100}],$ $[a_{000},$ $a_{100},$ $a_{110},$ $a_{010}],$ $[a_{000},$ $a_{001},$ $a_{011},$ $a_{010}],$ $[a_{111},$ $ a_{101},$ $a_{100},$ $a_{110}],$ $[a_{111},$ $a_{011},$ $a_{010},$ $a_{110}],$ $[a_{111},$ $a_{101},$ $a_{001},$ $a_{011}]\}$. The map $M$ is CS under the involution $I := a_{xyz}\mapsto a_{(1-x)(1-y(1-z)}$ where $x,y,z \in \{0,1\}$. We denote the faces $[a_{000},$ $a_{001},$ $a_{101},$ $a_{100}],$  $[a_{000},$ $a_{100},$ $a_{110},$ $a_{010}],$ $[a_{000},$ $a_{001},$ $a_{011},$ $ a_{010}],$ $[a_{111},$ $a_{101},$ $a_{100},$ $a_{110}],$ $[a_{111},$ $a_{011},$ $a_{010},$ $a_{110}],$ $[a_{111},$ $a_{101},$ $a_{001},$ $a_{011}]$ by $f_{1},$ $f_{2},$ $f_{3},$ $f_{4},$ $f_{6},$ $f_{5}$ respectively in dual $\overline{M}$ of $M$ which are vertices of $\overline{M}$. So, $\overline{M} := \{[f_{1},$ $f_{2},$ $f_{3}],$ $[ f_{1},$ $f_{3},$ $f_{5}],$ $[ f_{1},$ $f_{2},$ $f_{4}],$ $[f_{1},$ $f_{4},$ $f_{5}],$ $[ f_{2},$ $f_{3},$ $f_{6}],$ $[f_{2},$ $f_{4},$ $f_{6}],$ $[f_{3},$ $f_{5},$ $f_{6}],$ $[f_{4},$ $f_{5},$ $f_{6}]\}$. By the definition of CS, $[a_{000},$ $a_{001},$ $a_{101},$ $a_{100}]^{I}$ $=$ $[a_{111},$ $a_{011},$ $a_{010},$ $a_{110}]$, $[a_{000},$ $a_{100},$ $a_{110},$ $a_{010}]^{I}$ $= [a_{111},$ $a_{101},$ $a_{001},$ $a_{011}]$ and $[a_{000},$ $a_{001},$ $a_{011},$ $a_{010}]^{I} = [a_{111},$ $a_{101},$ $a_{100},$ $a_{110}]$. We define $\overline{I} := (f_{1},$ $f_{6})$ $(f_{2},$ $f_{5})$ $(f_{3},$ $ f_{4})$. The map $\overline{M}$ is centrally symmetric under $\overline{I}$ which is of type $\{3, 4\}$.

\section{Enumeration results for centrally symmetric $3$-manifolds}

We have modified the program MANIFOLD$_{-}$VT\cite{lutz1} and used in { Section 2} to enumerate CST surfaces. In this section we use the same modified program to enumerate CST { 3}-manifolds.
%We also use $\mathbb{Z}_{2} = <I>$ group action on the set $\{1, 2, \dots, 2m\}$ defined in $Section 2$.
It generates all possible { 2}- and { 3}-orbits. Let $F_{d}$ be a {\em d}-orbit for $d \in \{2, 3\}$. We ignore those { 3}-orbit for which $F_{3} \bigcap F_{3}^{I} \not= \emptyset$. Also we ignore those { 2}-orbit for which $F_{2} \bigcap F_{2}^{I} \not= \emptyset$. Hence we get all possible admissible { 2}- and { 3}-orbits. We check link of $ m$ vertices namely ${ 1, 2, \dots, m}$ which are used to define {\em I}. We also compute reduced homology groups of the objects using \cite{hom}. Hence we get all possible non isomorphic CST { 3}-manifolds. As a result for { $m = 6$} we have listed all possible { 3}-manifolds in Table \ref{table4}.

\bigskip

\noindent{\sc Proof of Theorem}\ref{thm5} The proof of theorem \ref{thm5} now follows from arguments in this section. Table \ref{table4} gives the list of centrally symmetric { 3}-manifolds on 12 vertices. The total number of objects is { 68}. By looking at homology groups we deduce that the objects are orientable . In this table we use used $a$ and $b$ to denote 10 and 11 respectively. We use {\em ijkl} to represent a 3-orbit in the list of orbits of order 2. where $i, j, k, l \in \{1, 2, \dots, 8, 9, a, b\}$.
\hfill$\Box$

\section{Acknowledgement}

Part of this work was done when first author was visiting Indian Institute of Science during summer of 2013. We gratefully acknowledges the support and guidance received from Prof. B. Datta. Work of second author is partially supported by SERB, DST grant No. SR/S4/MS:717/10.

\end{document}